\documentclass[11pt,twoside]{article}

\usepackage{lipsum} 

\usepackage{microtype} 

\usepackage{lmodern}
\usepackage[sans]{dsfont}
\usepackage[latin1]{inputenc}

\usepackage[hmarginratio=1:1,top=32mm,columnsep=20pt]{geometry} 
\usepackage{booktabs} 
\usepackage{hyperref} 

\usepackage{paralist} 

\usepackage{abstract} 

\usepackage{titlesec} 
\titleformat{\section}[block]{\large\scshape\centering}{\thesection.}{0.4em}{} 
\titleformat{\subsection}[block]{\large}{\thesubsection.}{0.4em}{} 

\let\OLDthebibliography\thebibliography
\renewcommand\thebibliography[1]{
  \OLDthebibliography{#1}
  \setlength{\parskip}{0pt}
  \setlength{\itemsep}{0pt plus 0.0ex}
    \advance\leftmargin20pt
}

\usepackage{fancyhdr} 
\usepackage{graphicx}
\usepackage{amsmath}
\usepackage{amssymb}
\usepackage{amsthm}
\usepackage{subfigure}
\usepackage{textcomp}
\usepackage{listings}
\usepackage{setspace}
\usepackage{graphicx}
\usepackage{epstopdf}
\usepackage{eurosym}
\usepackage{inputenc}
\usepackage{inputenc}
\usepackage{appendix}
\usepackage{tabularx}
\usepackage{array}
\usepackage{cite}
\usepackage{makeidx}
\usepackage{etoolbox}
\usepackage{mathtools}
\usepackage{empheq}
\usepackage{mathrsfs}
\usepackage{marginnote}

\DeclareGraphicsExtensions{.eps}
\DeclareGraphicsExtensions{.png}

\theoremstyle{plain}
\newtheorem{theorem}{Theorem}[section]
\newtheorem{lemma}[theorem]{Lemma}

\newtheorem{corol}[theorem]{Corollary}

\newtheorem{mydef}[theorem]{Definition}

\theoremstyle{definition}

\newtheorem{rem}[theorem]{Remark}

\hypersetup{citecolor=black}

\newenvironment{Proof}[1][Proof]
 {\proof[#1]\leftskip=0.0cm\rightskip=0.0cm}
 {\endproof}

\label{mathrefs}

\newcommand{\m}{\mbox{d}}

\newlength\tindent
\title{\vspace{-15mm}\fontsize{24pt}{10pt} \LARGE{\scshape{A nonlocal stochastic Cahn-Hilliard equation}}} 
\small{
\author{
\textsc{Federico Cornalba}
\\ \footnotesize{Department of Mathematical Sciences, University of Bath, Bath, BA2 7AY, United Kingdom}
\\ \vspace{0.7 pc}\footnotesize{e-mail: {\ttfamily F.Cornalba@bath.ac.uk}} \vspace{0.7 pc} 
\vspace{-5mm}
}}
\date{}

\begin{document}

\maketitle 
\thispagestyle{empty}

\pagestyle{fancy} 


\renewenvironment{abstract}
 {\small
  \begin{center}
  \bfseries \abstractname\vspace{-0.0 pc}\vspace{0pt}
  \end{center}
  \list{}{%
    \setlength{\leftmargin}{10mm}
    \setlength{\rightmargin}{\leftmargin}%
  }%
  \item\relax}
 {\endlist}

 \newenvironment{thanks}
 {
  \list{}{%
    \setlength{\leftmargin}{0mm}
    \setlength{\rightmargin}{\leftmargin}%
  }%
  \item\relax}
 {\endlist}

\begin{abstract}
\vspace{-0.5 pc}
We consider a stochastic extension of the nonlocal convective Cahn-Hilliard equation containing an additive Wiener process noise. We first introduce a suitable analytical setting and make some mathematical and physical assumptions. We then establish, in a variational context,  the existence of a \emph{weak statistical solution} for this problem. Finally we prove existence and uniqueness 
of a \emph{strong solution}.\vspace{0.5 pc}\vspace{0.5 pc}


\end{abstract}

\section{Introduction}

In the present work we study a stochastic version of the Cahn-Hilliard equation. This equation was initially used in \cite{CH1} to describe the evolution of the so called \emph{spinodal decomposition} phenomenon, a special type of phase separation process involving a binary metallic alloy. \\
The basic structure of the equation is 
\begin{equation}\label{eq:101}
\phi_t=-\mbox{div}\left\{-M(\phi)\nabla\mu(\phi)\right\}\quad\mbox{ in }D,
\end{equation}
where $D\subset\mathbb{R}^{d}$, $d\in\{2,3\}$, is a bounded domain with smooth boundary and $\phi$ is the unknown relative concentration of the metallic elements of the binary alloy. The term $M(\phi)$ is the \emph{mobility coefficient}. It determines the local speed of the mass flux in the alloy. The expression of $M(\phi)$ varies according to the accuracy of the physical modeling of this phenomenon. In this work we are not interested in discussing this term; we hence choose $M(\phi)$ to be the unitary constant. This is the easiest physically consistent choice. \\The term $\mu(\phi)$ is the \emph{chemical potential}. According to its nature, the Cahn-Hilliard equation assumes two different forms:
\begin{description}\itemsep1pt\parskip0pt\parsep0pt
\item (a) the \emph{local} form: if, for all $x\in D$, the quantity $\mu(\phi(x))$ can be computed by means of the values which $\phi$ takes in an arbitrarily small neighbourhood of $x$, then the Cahn-Hilliard equation is called local. A chemical potential possessing this property is
\begin{equation}\label{eq:920}
\mu=-\Delta\phi+F'(\phi),
\end{equation}
where $F:\mathbb{R}\rightarrow\mathbb{R}$ is a suitable function representing the \emph{density of potential energy}.
\item (b) the \emph{nonlocal} form: if $\mu$ does not satisfy the condition stated in (a), then the Cahn-Hilliard equation is called nonlocal. An example is the chemical potential
\begin{equation}\label{eq:921}
\mu=a\phi-J\ast\phi+F'(\phi),
\end{equation}
where $J:\mathbb{R}^d\rightarrow\mathbb{R}$ is a suitable kernel function and 
$$
(J\ast\phi)(x):=\int_{D}{J(x-y)\phi(y)\m y},\qquad a(x):=\int_{D}{J(x-y)\m y},\qquad\forall x\in D.
$$
\end{description}
\begin{rem}
The chemical potential \eqref{eq:920} was first introduced in \cite{CH1} by Cahn and Hilliard. The chemical potential \eqref{eq:921} was rigorously justified by Giacomin and Lebowitz in \cite{GL1,GL2}.
\end{rem}
\begin{rem}
The nonlocal Cahn-Hilliard equation is widely regarded as a better mathematical representation of the spinodal decomposition phenomenon than the local equation, the latter being a ``local approximation" of the nonlocal one. Nevertheless, the nonlocal equation is quite delicate to handle because of the presence of the convolution term. The mathematical aspects of this equation were studied, e.g., in \cite{BH1,BH2,GAJ,GZ,HIL,CFG,FG}.
\end{rem}
Equation \eqref{eq:101} gives an acceptable representation of the spinodal decomposition evolution. Nevertheless, it fails to take into account some ``dynamical" aspects of the alloy solute, such as vibrational, electronic and magnetic properties, as clarified in \cite[p.~217]{LHJSYK}. On the contrary, these aspects are important. In order to keep track of them, a statistical approach is required. See \cite{CO} for a detailed discussion.\\
The most direct way to include such aspects into \eqref{eq:101} is to add a Wiener process $w$ in it, thus obtaining the stochastic partial differential equation
\begin{equation}\label{eq:922}
\m\phi=-\mbox{div}\left\{-M(\phi)\nabla\mu(\phi)\right\}\m t+\sigma(\phi)\m w,
\end{equation}
where $\sigma$ is a suitable stochastic integrand. \\ \\
The \emph{local} version of \eqref{eq:922} has been studied in many papers. In \cite{EM}, existence of a weak solution and existence, uniqueness and measurability of a strong solution for the Cahn-Hilliard equation with additive constant stochastic noise (meaning that $\sigma$ is a multiple of the identity operator) are established by means of a variational approach. In \cite{DPD}, an analogous equation is examined with a pathwise method. Existence and uniqueness of a \emph{classical} solution is proved, along with existence and uniqueness of an invariant measure for the transition semigroup. This approach is completely different from the one followed in \cite{EM}. 
In \cite{CW}, an equation with noise having a nonlinear diffusion coefficient is considered. The author proves existence and uniqueness of a \emph{classical} solution, proves that this solution is differentiable in the sense of the Malliavian calculus, and, under some further assumptions, proves that the law of the solution is absolutely continuous with respect to Lebesgue measure.
In \cite{DZ}, a stochastic Cahn-Hilliard equation with reflection on a portion of the domain boundary $\Gamma$ and with the constraint of conservation of the space average is considered. Existence and uniqueness of a strong solution are shown for all continuous nonnegative initial conditions. Detailed information on the associated invariant measure and Dirichlet form is provided.
In \cite{AK}, the authors prove existence of solutions for a generalized stochastic Cahn-Hilliard equation with multiplicative white noise posed on bounded spatial convex domains, with piece-wise smooth boundary, and an additive time-dependent white noise term in the chemical potential. Existence is also derived for some non-convex cases when the boundary $\Gamma$ is smooth.\\ \\
In this paper we study a stochastic version of a specific \emph{nonlocal} Cahn-Hilliard equation of the form
$$
\m\phi=\left(-u\cdot\nabla\phi-\mbox{div}\left\{-M(\phi)\nabla\mu(\phi)\right\}\right)\m t+\sigma(\phi)\m w.
$$
More precisely, we have enriched \eqref{eq:922} with the \emph{convective} term $u\cdot\nabla\phi$, $u$ being a suitable \emph{velocity vector field}. Such modification is physically consistent and has been dealt with in the deterministic case, e.g., in \cite{CFG,FG}. In Sections \ref{sec:2} and \ref{sec:3} we set up the necessary notation and recall some useful results from probability and measure theory. In Section \ref{sec:4} we give the \emph{formal} definition of our problem. In Section \ref{sec:5} we introduce some mathematical and physical assumptions for our problem. In Section \ref{sec:40} we give the definition of \emph{weak statistical solution} for our problem and we prove the existence of such a solution, along with some further properties. In Section \ref{sec:41} we give the definition of \emph{strong} solution for our problem and we outline the strategy of the proof of existence and uniqueness of strong solutions. Finally, in Section \ref{sec:42} we provide the proofs outlined in the previous section and we thus achieve the existence and uniqueness 
of a strong solution.
\begin{rem}
The problem we will be examining is a stochastic version of \cite[p.~429, (1.9)/(1.10)/(1.13)]{CFG}. In addition, the structure of Sections \ref{sec:40} and \ref{sec:42} is partially inspired by \cite{EM,FV}.
\end{rem}

\section{Basic Notation}\label{sec:2}
We define a few mathematical objects whose notation will be kept throughout the entire paper.
\begin{itemize}\itemsep1pt\parskip0pt\parsep0pt
\item The bounded spatial domain where our specific Cahn-Hilliard equation lives is denoted by $D$, where $D\subset\mathbb{R}^d$, $d\in\{2,3\}$. Its boundary is indicated as $\Gamma$. The time interval in which we make all our considerations is $[0,T]$, for an arbitrary given $T>0$.
\item The symbol $L^p$ (where $p\in[1;+\infty]$) denotes both $L^p(D)$ and $\left[L^p(D)\right]^d$. The symbol $H^s$ (where $s\geq 0$) denotes both $H^s(D)$ and $\left[H^s(D)\right]^d$. For $s\geq 0$, we denote by $H^{-s}$ the dual space of $H^s$. For a Banach space $W$, we denote by $L^p([0,T];W)$ (where $p\in[1,+\infty]$) the $L^p$-space of $W$-valued, Bochner-integrable functions on $[0,T]$. We denote by $\mathcal{C}([0,T];W)$ (resp. $\mathcal{C}^{\beta}([0,T];W)$) the set of $W$-valued continuous (resp. $\beta$-H\"older continuous) functions on $[0,T]$, for $\beta\in(0,1)$. Moreover, for $k\in\mathbb{N}$ and $p\in[1,+\infty]$, we set $W^{k,p}:=\left\{v\in L^p:D^{\alpha}u\in L^p\mbox{ for each multi-index }\alpha\mbox{ such that }|\alpha|\leq k\right\}$. The spaces $W^{k,p}(\mathbb{R}^d)$ have an analogous definition with $\mathbb{R}^d$ replacing $D$.
\item We set
$$
H:=L^2,\qquad U:=H^1,\qquad V:=\left\{u\in H^2:\frac{\partial u}{\partial \nu}=0\mbox{ on }\Gamma\right\}\!.
$$
We denote by $(\cdot,\cdot)$ and $\|\cdot\|$ the inner product and norm in $H$, respectively.
\item For a Banach space $W$ which is not $H$, we denote its dual space by $W'$, its duality pairing by $\langle \cdot,\cdot\rangle_{W',W}$, and its norm by $\|\cdot\|_{W}$. If $W$ is also a Hilbert space, we denote its inner product by $(\cdot,\cdot)_{W}$.
\item For a Banach space $Z$ and a subset $W\subset Z$, we indicate the Borel $\sigma$-algebra of $W$ with respect to the topology of $Z$ as $\mathcal{B}_{Z}(W)$, namely $\mathcal{B}_{Z}(W):=\mathcal{B}(Z)\cap W$, where $\mathcal{B}(Z)$ denotes the standard Borel $\sigma$-algebra of $Z$. The previous definitions are consistent if $W\in\mathcal{B}(Z)$.
\item $(\Omega,\mathcal{F},\mathbf{m})$ denotes the main probability space. The symbol $\mathbb{E}$ stands for the \emph{expected value} operator. 
We use bold characters to indicate measures on measurable spaces. We use the bold characters along with a hat ($\hat{\cdot}$) to indicate the \emph{characteristic functionals} of measures. 
\item We use the concise notations a.e. and a.s. to say \emph{almost everywhere} and \emph{almost surely}. 
We will mainly use the notation a.e. when referring to measures on spatial or time domains and the notation a.s. when referring to probability measures on a probability space.
\end{itemize}
The definition of characteristic functionals and of some more specialized analytical and probabilistic tools will be specified in the following section.

\section{preliminary tools}\label{sec:3}
We here recall some useful definitions and results related to probability and measure theory.
\begin{mydef}[Characteristic functional]\label{def:902}
Given a separable Banach space $Y$ and a probability measure $\nu$ defined on $(Y,\mathcal{B}(Y))$, the {\bfseries characteristic functional} of $\nu$ is a $\mathbb{C}$-valued functional with domain $Y'$, defined as
$$
\hat{\nu}(f):=\int_{Y}{\exp\left\{i\langle f , y\rangle_{Y',Y}\right\}\nu(\mbox{\emph{d}} y)},\qquad\forall f\in Y'.
$$
If $X$ is a $Y$-valued random variable, the {\bfseries characteristic functional} of $X$ is the characteristic functional of the law of $X$, namely
$$
\hat{\mathcal{L}}_X(f):=\int_{Y}{{\exp\left\{i\langle f , y\rangle_{Y',Y}\right\}\mathcal{L}_{X}(\mbox{\emph{d}} y)}}=\mathbb{E}\left[\exp\{i\langle f,X\rangle_{Y',Y}\}\right],\qquad \forall f\in Y',
$$
where $\mathcal{L}_X$ denotes the law of $X$. If $Y$ is a separable Hilbert space, we adapt the previous expressions by replacing the duality $\langle\cdot,\cdot\rangle_{Y',Y}$ with the $Y$-inner product $(\cdot,\cdot)_{Y}$ and by considering $f\in Y$ instead of $f\in Y'$.
\end{mydef}
\begin{mydef}[Pushforward measures]\label{def:901}
Let $(X_1,\mathcal{F}_1)$, $(X_2,\mathcal{F}_2)$ be two measurable spaces. Let $f:X_1\rightarrow X_2$ be a measurable function and let $\nu$ be a probability measure on $(X_1,\mathcal{F}_1)$. The {\bfseries pushforward measure} of $\nu$ associated to $f$ is the probability measure $\nu_{\ast}$ on $(X_2,\mathcal{F}_2)$ defined as
$$
\nu_{\ast}(C):=\nu\left(f^{-1}(C)\right)\!,\quad\forall C\in\mathcal{F}_2.
$$
\end{mydef}
\begin{lemma}\label{lem:903}
With the notation adopted in Definition \eqref{def:901}, let $g:X_2\rightarrow\mathbb{C}$ be a measurable function. Then $g$ is integrable on $X_2$ with respect to $\nu_{\ast}$ if and only if $g\circ f$ is integrable on $X_1$ with respect to $\nu$. In this case, the integrals coincide, i.e.
\begin{equation}\label{eq:905}
\int_{X_2}{g\,\mbox{\emph{d}}\nu_{\ast}}=\int_{X_1}{(g\circ f)\,\mbox{\emph{d}}\nu}.
\end{equation}
\end{lemma}
\begin{theorem}[It\"o formula]\label{th:904}
Let $Y$ be a Hilbert space. Let $(\Omega,\mathcal{F},\mathbf{m},\{\mathcal{F}_t\}_{t\geq 0},\{w_t\}_{t\geq 0})$ be a $Y$-valued $Q$-Wiener process (in the sense of \emph{\cite[p.~13, Definition 2.1.9.]{PR}}). Let $\varphi$ be a $Y$-valued, $[0,T]$-Bochner integrable, predictable process. Let $\Phi$ be a $L^0_{2}$-valued process stochastically integrable in $[0,T]$. Let $X(0)$ be a $\mathcal{F}_0$-measurable $Y$-valued random variable. Consider the (well defined) process 
$$
X(t):=X(0)+\int_{0}^{t}{\varphi(s)\mbox{\emph{d}} s}+\int_{0}^{t}{\Phi(s)\mbox{\emph{d}} w(s)}.
$$
Assume that a function $F:[0,T]\times Y\rightarrow\mathbb{R}$ and its partial Fr\'echet derivatives $F_t,F_y,F_{yy}$ are uniformly continuous on bounded subsets of $[0,T]\times Y$. Then the following {\bfseries It\"o formula} holds a.s. for all $t\in[0,T]$
\begin{eqnarray*}
F(t,X(t)) & = & F(0,X(0))+\int_{0}^{t}{\left(F_y(s,X(s)),\Phi(s)\mbox{\emph{d}} w(s)\right)}
+\int_{0}^{t}{\left\{F_t(s,X(s))+\left(F_y(s,X(s)),\varphi(s)\right)\right\}\mbox{\emph{d}} s}\\& + &\int_{0}^{t}{\left\{\frac{1}{2}\mbox{\emph{tr}}\left[F_{yy}(s,X(s))\left(\Phi(s)Q^{1/2}\right)\left(\Phi(s)Q^{1/2}\right)^{\ast}\right]\right\}\mbox{\emph{d}} s},
\end{eqnarray*}
where $Q$ is the covariance operator of the ($\,Y$-valued) Wiener process $w$. 
\end{theorem}
We also refer the reader to  \cite[Section 2.1.]{PR} and \cite[p.~36, Lemma 2.4.2.]{PR} for a detailed discussion of infinite-dimensional Wiener processes and to \cite[Paragraph 4.5]{DPZ} for a detailed discussion of the It\"o formula.

\section{Abstract formulation of the problem}\label{sec:4}
We \emph{formally} study the stochastic partial differential equation
\begin{subequations}
\label{eq:303}
\begin{empheq}[left={}\empheqlbrace]{align}
  & \,\,\m\phi=(-u\cdot\nabla\phi+\Delta\mu)\m t+\m w, \label{eq:304}\\
  & \,\,\mu=a\phi-J\ast \phi+F'(\phi),\label{eq:305}\\
  & \,\,\frac{\partial\mu}{\partial \nu}=0\quad\mbox{ on }\Gamma\times(0,T),\label{eq:337}\\
  & \,\,\phi(0)=\phi_0\in U,\nonumber
\end{empheq}
\end{subequations}
where $w$, $u$, $J$, $a$, $F$, $\phi_0$ are mathematical objects whose nature will be specified in the following section. The symbol $\ast$ in \eqref{eq:305} denotes the convolution operator over $D$, namely
$$
(J\ast\phi)(x):=\int_{D}{J(x-y)\phi(y)\m y},\quad\forall x\in D.
$$
Problem \eqref{eq:303} is a stochastic extension of \cite[p.~429, (1.9)/(1.10)/(1.13)]{CFG}.

\section{Assumptions}\label{sec:5}
We now give a precise meaning to the objects appearing in Problem \eqref{eq:303}. More precisely, we work under the following mathematical and physical assumptions.
\begin{description}
\item (i)\label{pag:10} $u$ is a given velocity field satisfying $u\in L^{\infty}([0,T]\times D)$, $\mbox{div}(u)=0$ in $D$, $u=0$ on $\Gamma$. 
\item (ii) $w=w(t)$, $t\in[0,T]$, is a $H$-valued $Q$-Wiener process in the sense of \cite[p.~13, Definition 2.1.9.]{PR}. Here \label{pag:30} $Q:H\rightarrow H$ is a continuous, symmetric, positive definite, finite trace linear operator.
\item (iii) $J$ is a kernel function 
satisfying the following properties:
$$
J\in W^{1,1}(\mathbb{R}^d),\qquad J(x)=J(-x),\qquad a(x):=\int_{D}{J(x-y)\m y}\geq 0,\quad\mbox{for a.e. }x\in D.
$$ 
\item (iv)\label{pag:100} We choose the density of potential energy 
$F$ to be 
$$
F(s)=\frac{s^4}{4}-\frac{s^2}{2},
$$
so that $F'(s)=s^3-s$ and $F\in\mathcal{C}^{2,1}_{loc}(\mathbb{R})$. We assume that there exists $c_0>0$ such that
\begin{equation}\label{eq:353}
F''(s)+a(x)\geq c_0,\quad\forall s\in\mathbb{R},\quad\mbox{a.e. }x\in D.
\end{equation}
The properties we have just listed for $F$ and $J$ are exactly hypotheses (H1)-(H3) in \cite[p.~431]{CFG}. 
\item (v) $H$ has an orthonormal basis $\{e_i\}_{i\in\mathbb{N}}$ consisting of the eigenvectors of the operator $A:\mathcal{D}(A)=V:v\mapsto(-\Delta+I)v$ with associated eigenvalues $\{\mu_i\}_{i\in\mathbb{N}}$, i.e., $Ae_i=\mu_ie_i$, $\forall i\in\mathbb{N}$. We recall that $\mu_i\geq 1,\forall i\in\mathbb{N}$, and that $\mu_i\rightarrow +\infty$. The space $V$ is endowed with the norm
$$
\|v\|_{V}:=\left(\|v\|^2+\|\Delta v\|^2\right)^{1/2}\!\!\!\!\!\!\!,\qquad \forall v\in V.
$$
Due to the regularity of $D$, this norm is equivalent to the standard $H^2$-norm. It follows that $\{e_i\}_{i\in\mathbb{N}}$ is an orthogonal basis in $V$. The family $\{e_i\}_{i\in\mathbb{N}}$ is also an orthogonal basis in $U$. In addition, we indentify $H$ with its dual space by means of the Riesz isomorphism and use the continuous injections
$$
V\hookrightarrow H^{\varepsilon}\hookrightarrow H\equiv H'\hookrightarrow H^{-\varepsilon}\hookrightarrow V',
$$
for a chosen parameter $\varepsilon\in(0,1/4)$.
\item (vi) $D\subset\mathbb{R}^{d}$, $d\in\{2,3\}$, is regular enough to apply the Rellich-Kondrachov Theorem 
and \cite[p.~1285, Theorem 1]{GR}. As consequences we have that $U$ is compactly embedded in $L^4$, that $\{e_i\}_{i\in\mathbb{N}}\subset\mathcal{C}^{\infty}(\overline{D})$, and that
$$
\|e_i\|_{L^{\infty}}\leq C(D)(\mu_i-1)^{(d-1)/4}.
$$
$D$ is also regular enough that $H$ is compactly embedded in the interpolation space $[H,V']_{\varepsilon/2}=H^{-\varepsilon}$. For details on this embedding, the reader is referred to \cite[pages 99-103]{LM}.
\item (vii) $\{e_i\}_{i\in\mathbb{N}}$ are eigenvectors for $Q$ as well, namely there is sequence of nonnegative real numbers $\{\vartheta_i\}_{i\in\mathbb{N}}$ such that $Qe_i=\vartheta_ie_i$ for each $i\in\mathbb{N}$. In addition we require that
\begin{equation}\label{eq:329}
K(Q):=\sum_{i=1}^{\infty}{(\mu_i-1)^{(d-1)/2}\vartheta_i}<+\infty.
\end{equation}
\item (viii)\label{pag:200} $\phi_0$ is a $U$-valued random variable which is independent of $w$. 
\end{description}
\begin{rem}
Due to the nature of $w$, equation \eqref{eq:304} is a stochastic infinite-dimensional differential equation. It can be \emph{formally} interpreted, with the due careful analogies and generalizations, as in \cite[p.~73, Definition 4.2.1.]{PR} or \cite[Chapter 7]{DPZ}.
\end{rem}
\begin{rem}\label{rem:12}
Condition \eqref{eq:329} is stronger than requiring $Q$ to have finite trace. However, it enables us not to make any harmful assumptions on the geometry of $\Gamma$. The geometry of $\Gamma$, in fact, affects many interpolation results in Sobolev spaces.\\
We could have replaced condition \eqref{eq:329} with a condition of uniform boundedness of the family $\{e_i\}_{i\in\mathbb{N}}$ in $L^{\infty}$; by so doing, however, we would have been forced to require additional conditions on the geometry of $\Gamma$ and, consequentially, we would have had to check the validity of some interpolation results. The latter approach is very tough and hence inadvisable.
\end{rem}
\begin{rem}
We could have chosen many other polynomial growths for the energy $F$ without disrupting its physical meaning. See \cite{DPD} for instance. Nevertheless, the fourth degree growth given in assumption (iv) plays a crucial role in some of the forthcoming mathematical computations. See for instance Remark \eqref{rem:199} below.
\end{rem}

\section{Existence of a weak statistical solution}\label{sec:40}
In this section we prove the existence of a \emph{weak statistical solution} to Problem \eqref{eq:303}, as defined in Definition \eqref{def:2} below.\\
We first introduce some function spaces for a given time $T>0$. Let
$$
\mathscr{U}:=L^2([0,T];U)\cap L^{\infty}([0,T];H)\cap\mathcal{C}^{2/5}([0,T];V')\cap L^4([0,T];L^4),
$$
$$
\mathscr{Z}:=L^{p'}([0,T];L^4)\cap\mathcal{C}([0,T];H^{-\varepsilon}),\quad p'\in(3,4),\quad\varepsilon\in(0,1/4),
\vspace{0.4 pc}
$$
where $p'$ is a chosen parameter and $\varepsilon$ has been introduced in Section \ref{sec:5}, hypothesis (v).\\
Because of the compatible nature of the Banach spaces appearing in the definitions of $\mathscr{U}$ and $\mathscr{Z}$, we can define their norms in the following way
$$
\|v\|_{\mathscr{U}}:=\|v\|_{L^2([0,T];U)}+\|v\|_{L^{\infty}([0,T];H)}+\|v\|_{\mathcal{C}^{2/5}([0,T];V')}+\|v\|_{L^4([0,T];L^4)},\,\,\,\forall v\in\mathscr{U}\!\!,
$$
$$
\|v\|_{\mathscr{Z}}:=\|v\|_{L^{p'}([0,T];L^4)}+\|v\|_{\mathcal{C}([0,T];H^{-\varepsilon})},\,\,\,\forall v\in\mathscr{Z}.
\vspace{0.4 pc}
$$
\vspace{0.2 pc}
\hspace{-0.2 pc}We state and prove a preliminary result.
\begin{theorem}\label{th:20}
$\mathscr{U}$ is compactly embedded in $\mathscr{Z}$.
\end{theorem}
\begin{Proof}
Let $\mathcal{S}$ be a bounded set in $\mathscr{U}$. If we apply \cite[p.~86, Theorem 6]{S} with $X=U$, $B=L^4$, $Y=V'$, $q=4$, $p=p'$, we deduce that $\mathcal{S}$ is relatively compact in $L^{p'}([0,T];L^4)$. If we apply \cite[p.~84, Theorem 5]{S} with $X=H$, $B=[H,V']_{\varepsilon/2}=H^{-\varepsilon}$ ($\varepsilon\in(0,1/4)$), $Y=V'$, $p=\infty$, we deduce that $\mathcal{S}$ is relatively compact in $\mathcal{C}([0,T];H^{-\varepsilon})$, 
hence the conclusion. 
\end{Proof}
For each $m\in\mathbb{N}$, we define the spaces
$$
V_m:=\mbox{span}\{e_1,\cdots,e_m\},
$$
\vspace{-1 pc}
$$
\mathscr{V}_m := \left\{v\in L^2([0,T];V)\cap L^{q'}([0,T];W^{2,4})\cap H^2([0,T];H)\!:\,
v(t)\in V_m\mbox{ for a.e. }t\in[0,T],\,\,v(T)=0\right\}\!,
$$
where $H^2([0,T];H)$ is the standard Sobolev time-dependent space and
\begin{equation}\label{eq:318}
\frac{2}{3}+\frac{1}{p'}+\frac{1}{q'}=1.
\end{equation}
\begin{rem}\label{rem:21}
We recall that $\cup_{r=1}^{\infty}{V_r}$ is dense in $H^{\varepsilon}\!$. As a consequence, $H^{\varepsilon}$ is separable.
\end{rem}
We now define a proper \emph{test function} space for our problem as 
\vspace{-0.5 pc}
$$
\mathscr{V}:=\mbox{ completion of }\bigcup_{m=1}^{\infty}{\mathscr{V}_m}\mbox{ with respect to the norm }\|\cdot\|_{\mathscr{V}},\mbox{ where}
$$
\vspace{-1.4 pc}
$$
\|v\|^2_{\mathscr{V}}:=\|v\|^2_{L^2([0,T];V)}+\|v\|^2_{L^{q'}([0,T];W^{2,4})}+\|v\|^2_{H^2([0,T];H)}.
$$
We are now able to define an important object related to $w$.
\begin{mydef}
We define the {\bfseries white noise} $\partial w/\partial t$ as the distributional time derivative of the Wiener process $w$. Namely, $\partial w/\partial t$ is a $\mathscr{V}'$-valued random variable such that
\begin{equation}\label{eq:360}
\left\langle \frac{\partial w}{\partial t}(\omega), v\right\rangle_{\!\!\mathscr{V}',\mathscr{V}}:=-\left(w(\omega),\frac{\partial v}{\partial t}\right)_{\!\!L^2([0,T];H)}\!\!,\!\!\qquad \forall v\in\mathscr{V}.
\end{equation}
\end{mydef}
\begin{rem}\label{rem:20}
The operator 
$$
\mathcal{A}:L^2([0,T];H)\rightarrow \mathscr{V}':w\mapsto \frac{\partial w}{\partial t},
$$
with the time derivative defined as in \eqref{eq:360}, is continuous, hence measurable. Using \cite[p.~65, Theorem 10.1., item (c)]{JP} we deduce that $\phi_0$ and $\partial w/\partial t$ are independent. In addition, the mapping
$$
\omega\mapsto\left\{\phi_0(\omega),\frac{\partial w}{\partial t}(\omega)\right\}
$$
is a  random variable from $(\Omega,\mathcal{F})$ to $(H^{-\varepsilon}\times \mathscr{V}',\mathcal{B}(H^{-\varepsilon}\times \mathscr{V}'))$.
\end{rem}
\begin{rem}\label{rem:23}
The test function space $\mathscr{V}$ is reflexive and separable. To see this, let us consider the space
$$
\mathscr{W}:=L^2([0,T];V)\cap L^{q'}([0,T];W^{2,4})\cap H^2([0,T];H).
$$
The space $\mathscr{W}$ is clearly reflexive and separable, being the intersection of reflexive and separable Banach spaces. Since $\mathscr{V}$ is a closed subspace of $\mathscr{W}$, it is reflexive and separable.
\end{rem}
We now specify what we mean by \emph{weak statistical solution}.
\begin{mydef}\label{def:2}
A {\bfseries weak statistical} solution (or simply a {\bfseries weak} solution) to Problem \eqref{eq:303} is a probability measure $\mathbf{P}$ (concentrated) on $\mathcal{B}_{\mathscr{Z}}(\mathscr{U})$ which, for every $\xi\in H^{\varepsilon}$ and $v\in\mathscr{V}\!$, satisfies
\begin{equation}\label{eq:323}
\int_{\mathscr{Z}}{\exp\left\{ i\langle \phi(0),\xi\rangle_{H^{-\varepsilon},H^{\varepsilon}}+iC(\phi,v)\right\}\mathbf{P}(\emph{d} \phi)}=\hat{\Xi}(\xi)\hat{\mathbf{N}}(v),
\end{equation}
where 
\begin{eqnarray*}
C(\phi,v) :=  -(\phi_0,v(0))-\int_{0}^{T}{\!\!\!\int_{D}{\phi u\cdot\nabla v}}-\int_{0}^{T}{\!\!\!\int_{D}{(a\phi+\phi^3-\phi-J\ast\phi)\Delta v}}-\int_{0}^{T}{\!\!\left(\phi,\frac{\partial v}{\partial t}\right)}.
\end{eqnarray*}
Here $\Xi$ indicates
the distribution of the random variable $\phi_0$ on $H$, and $\hat{\mathbf{N}}$ is a functional on $\mathscr{V}$ defined by
$$
\hat{\mathbf{N}}(v):=\hat{\mathbf{W}}\left(-\frac{\partial v}{\partial t}\right)\!,\quad\forall v\in\mathscr{V},
$$
where $\mathbf{W}$ is the distribution of $w$. 
\end{mydef}
The above notation $\hat{\mathbf{N}}$ is justified by the easily observed fact that $\hat{\mathbf{N}}$ is the characteristic functional of the white noise $\partial w/\partial t$. Definition \eqref{def:2} is analogous to \cite[p.~1181, Definition 5.1.]{EM}.
\begin{rem}
A weak solution is defined (more precisely, concentrated) on the $\sigma$-algebra $\mathcal{B}_{\mathscr{Z}}(\mathscr{U})$, since $\mathscr{U}$ is not separable. For a similar discussion, see \cite[p.~1181]{EM}. This solution may not be unique. In addition, we warn the reader not to confuse $\phi(0)$ and $\phi_0$ in Definition \eqref{def:2}.
\end{rem}
\begin{rem}\label{rem:16}
Here we give a \emph{formal} 
derivation of Definition \eqref{def:2}. If $\phi$ is a solution of the infinite dimensional stochastic equation \eqref{eq:304}-\eqref{eq:305} (e.g., in the sense of \cite[p.~73, Definition 4.2.1.]{PR}), then we may compute, by means of the It\"o formula, the stochastic differential equation which is satisfied by the real-valued process $B(\phi):=(\phi,v)$, for each $v\in\cup_{m=1}^{\infty}\mathscr{V}_m$. We obtain
$$
0= (\phi_0,v(0))-\int_{0}^{T}{\!\!\int_{D}{vu\cdot\nabla\phi}}+\int_{0}^{T}{\!\!\int_{D}{(\Delta\mu)v}}+\int_{0}^{T}{\left(\phi,\frac{\partial v}{\partial t}\right)}-\int_{0}^{T}{\left(\frac{\partial v}{\partial t},w\right)}.
$$
If we recall the definition of white noise and condition \eqref{eq:337}, and we observe that $(v,u\cdot\nabla\phi)=-(\phi,u\cdot\nabla v)$ thanks to the properties of $u$ listed in assumption (i), Section \ref{sec:5}, we derive
\begin{equation}\label{eq:338}
C(\phi,v)=\left\langle \frac{\partial w}{\partial t},v\right\rangle_{\mathscr{V}',\mathscr{V}}\!\!\!\!\!\!\!\!\!\!\!\!.
\end{equation}
Relation \eqref{eq:338} can be extended by a density argument (with respect to the norm $\|\cdot\|_{\mathscr{V}}$) to all $v\in\mathscr{V}$, provided that $\phi$ is sufficiently regular.\\
If we add the term $\langle \phi(0),\xi\rangle_{H^{-\varepsilon},H^{\varepsilon}}$ to both sides of \eqref{eq:338}, we multiply by $i$, apply the exponential function, take the expected value and use the independence of $w$ and $\phi_0$, we see that \eqref{eq:323} holds with $\mathbf{P}$ being the distribution of $\phi$. We may hence generalize \eqref{eq:323} omitting the requirement that $\mathbf{P}$ is the distribution of a given process, thus obtaining Definition \eqref{def:2}. Equality \eqref{eq:338} also plays an important role in the forthcoming Definition \eqref{def:4}, aside from the fact that it justifies the expression of $C(\phi,v)$.\\ We stress the fact that the computations made in this remark are, for the time being, formal, and they will later be suitably justified in a rigorous context. 
\end{rem}
\begin{rem}\label{rem:10}
It can be seen that the real-valued functional $\phi\mapsto C(\phi,v)$ is continuous on $\mathscr{Z}$ for each fixed $v\in\mathscr{V}$.
To show this, let $\phi_n$ be a sequence such that $\phi_n\rightarrow \phi$ in $\mathscr{Z}$. Because of the convergence in $L^{p'}([0,T];L^4)$, it is straightforward to deduce the suitable convergence of all the elements appearing in $C(\phi_n,v)$ which are linear in $\phi_n$. We only have to treat the nonlinearity $\phi^3$ with a little bit of care. In fact, recalling \eqref{eq:318} and using the H\"older inequality first in space and then in time, we obtain
\begin{eqnarray}\label{eq:352}
\left|\int_{0}^{T}{\!\!\!\!\int_{D}{(\phi^3_n-\phi^3)}\Delta v}\right| & \leq & \frac{3}{2}\!\int_{0}^{T}{\!\!\!\!\int_{D}{\!|\phi_n-\phi|(\phi_n^2+\phi^2)|\Delta v|}}
 \leq  \frac{3\sqrt{2}}{2}\!\int_{0}^{T}{\!\|\phi_n-\phi\|_{L^4}(\|\phi_n\|^2_{L^4}+\|\phi\|^2_{L^4})\|\Delta v\|_{L^4}}\nonumber\\
& \leq & 3\|\phi_n-\phi\|_{L^{p'}([0,T];L^4)}\left[\|\phi_n\|^2_{L^3([0,T];L^4)}+\|\phi\|^2_{L^3([0,T];L^4)}\right]\|\Delta v\|_{L^{q'}([0,T];L^4)}\nonumber\\
& \leq & 3\|\phi_n-\phi\|_{L^{p'}([0,T];L^4)}\left[\|\phi_n\|^2_{L^3([0,T];L^4)}+\|\phi\|^2_{L^3([0,T];L^4)}\right]\|v\|_{\mathscr{V}}.
\end{eqnarray}
We deduce
$$\left|\int_{0}^{T}{\!\!\!\int_{D}{(\phi^3_n-\phi^3)}\Delta v}\right|\rightarrow 0
$$
as $n\rightarrow+\infty$, hence the conclusion. 
Inequality \eqref{eq:352} also allows us to define $C(\phi)\in\mathscr{V}'$ as 
$$
\langle C(\phi),v\rangle_{\mathscr{V}',\mathscr{V}}:=C(\phi,v)
$$
for every $\phi\in\mathscr{Z}$.
\end{rem}
\begin{rem}\label{rem:13}
It is straightforward to deduce two facts from Remark \eqref{rem:10}.
\begin{description}
\item (a) For each $\xi\in H^{\varepsilon}$ and each $v\in\mathscr{V}$, the functional 
$$
\phi\mapsto\exp\left\{ i\langle \phi(0),\xi\rangle_{H^{-\varepsilon},H^{\varepsilon}}+iC(\phi,v)\right\}
$$
is continuous on $\mathscr{Z}$.
\item (b) The mapping 
$$
\mathscr{D}:\mathscr{Z}\rightarrow H^{-\varepsilon}\times\mathscr{V}':\phi\mapsto \{\phi(0),C(\phi)\}
$$
is continuous as a mapping from $(\mathscr{Z},\mathcal{B}(\mathscr{Z}))$ to $(H^{-\varepsilon}\times\mathscr{V}',\mathcal{B}(H^{-\varepsilon}\times\mathscr{V}'))$. This is consequence of the nature of $\|\cdot\|_{\mathscr{V}}$ and of \eqref{eq:352}. The mapping $\mathscr{D}$ \emph{simultaneously} takes into account the initial condition $\phi(0)$ and the ``deterministic body'' $C(\phi)$ associated with our Cahn-Hilliard equation \eqref{eq:304}. This compound nature of $\mathscr{D}$ will play a crucial role in Section \ref{sec:42} for our conclusive and most important result.
\end{description}
\end{rem}
In the forthcoming sections, we will need the following fact.
\begin{lemma}\label{lem:70}
Equality \eqref{eq:323} implies that 
\begin{equation}\label{eq:330}
\mathbf{P}(\mathscr{D}^{-1}(C))=(\Xi\times\mathbf{N})(C),\quad\forall C\in\mathcal{B}(H^{-\varepsilon}\times\mathscr{V}').
\end{equation}
The left hand side of equality \eqref{eq:330} is well defined since $\mathscr{D}$ is continuous, as stated in Remark \eqref{rem:13}.
\end{lemma}
\begin{Proof}
Let $\xi\in H^{\varepsilon}$, $v\in\mathscr{V}$ be arbitrarily fixed. In accordance with the notation introduced in Definitions \eqref{def:901} and \eqref{def:2}, we set $X_1:=\mathscr{Z}$, $X_2:=H^{-\varepsilon}\times\mathscr{V}'$, $f:=\mathscr{D}$, $\nu:=\mathbf{P}$, and
$$
g:X_2\rightarrow\mathbb{C}:\left\{x,y\right\}\mapsto\exp\left\{i\langle x,\xi \rangle_{H^{-\varepsilon},H^{\varepsilon}}+i\langle y, v\rangle_{\mathscr{V}',\mathscr{V}}\right\}\!.
$$
If we use \eqref{eq:905}, and the definition of weak statistical solution, 
we deduce
\begin{eqnarray}\label{eq:365}
& &\!\!\!\!\!\int_{\mathscr{Z}}{\exp\left\{ i\langle \phi(0),\xi\rangle_{H^{-\varepsilon},H^{\varepsilon}}+iC(\phi,v)\right\}\mathbf{P}(\m\phi)}=\int_{X_1}{(g\circ f)\,\m\nu}=\int_{X_2}{g\,\m\nu_{\ast}}\nonumber\\
&=&\int_{H^{-\varepsilon}\times\mathscr{V}'}{\!\!\!\!\!\!\!\!\exp\left\{i\langle x,\xi\rangle_{H^{-\varepsilon},H^{\varepsilon}}+i\langle y,v \rangle_{\mathscr{V}',\mathscr{V}}\right\}\mathbf{P}_{\ast}(\m\!\left\{x,y\right\})}=\hat{\Xi}(\xi)\hat{\mathbf{N}}(v)\nonumber\\
&=&\int_{H^{-\varepsilon}\times\mathscr{V}'}{\!\!\!\!\!\!\!\!\exp\left\{i\langle x,\xi\rangle_{H^{-\varepsilon},H^{\varepsilon}}+i\langle y,v \rangle_{\mathscr{V}',\mathscr{V}}\right\}\left(\Xi\times\mathbf{N}\right)(\m\!\left\{x,y\right\})}.
\end{eqnarray}
We have already noticed in Remarks \eqref{rem:21}, \eqref{rem:23} that $H^{\varepsilon}$ and $\mathscr{V}$ are reflexive and separable, thus $H^{\varepsilon}\times\mathscr{V}$ is reflexive and separable. Since $\xi$ and $v$ are arbitrarily chosen in $H^{\varepsilon}$ and $\mathscr{V}$, the reflexivity of $H^{\varepsilon}$ and $\mathscr{V}$, and relation \eqref{eq:365} imply that 
\begin{equation}\label{eq:366}
\int_{H^{-\varepsilon}\times\mathscr{V}'}{\!\!\!\!\!\exp\left\{iL(\left\{x,y\right\})\right\}\mathbf{P}_{\ast}(\m\!\left\{x,y\right\})}=\int_{H^{-\varepsilon}\times\mathscr{V}'}{\!\!\!\!\!\exp\left\{iL(\{x,y\})\right\}(\Xi\times\mathbf{N})(\m\!\left\{x,y\right\})}
\end{equation}
for each $L\in\left(H^{-\varepsilon}\times\mathscr{V}'\right)'$.
Because $H^{\varepsilon}\times\mathscr{V}$ is reflexive and separable, it follows that $\left(H^{\varepsilon}\times\mathscr{V}\right)'$ is separable. Since there is an isometric isomorphism between $\left(H^{\varepsilon}\times\mathscr{V}\right)'$ and $H^{-\varepsilon}\times\mathscr{V}'$, we have that $H^{-\varepsilon}\times\mathscr{V}'$ is separable. We can hence apply \cite[p.~28, Proposition 4.15.]{AG} and deduce that $\mathbf{P}_{\ast}\equiv\left(\Xi\times\mathbf{N}\right)$, i.e. \eqref{eq:330}.
\end{Proof}
We can now state and prove the main theorem of this section.
\begin{theorem}\label{th:19}
Let $d\leq 3$. Let $w$ be a $H$-valued $Q$-Wiener process and let $w$, $u$, $Q$, $J$, $F$, $\phi_0$ and $\{e_j\}_{j\in\mathbb{N}}$ satisfy the properties \emph{(i)-(viii)} listed in Section \ref{sec:5}. Let $\phi_0$ be a $U$-valued random variable satisfying
$$
\mathbb{E}\left[\|\phi_0\|_{U}^2+\int_{D}{\frac{\phi^4_0}{4}}-\int_{D}{\frac{\phi^2_0}{2}}\right]<+\infty.
$$
Then Problem \eqref{eq:303} admits a weak statistical solution in the sense of Definition \eqref{def:2}.
\end{theorem}
\begin{Proof}[Proof of Theorem \eqref{th:19}.]
\emph{Step 1: Galerkin Approximation of Problem \eqref{eq:303}}. For each $m\in\mathbb{N}$, we denote by $\pi_m$ the $H$-orthogonal projection operator on $V_m$. More precisely, we use the extended operator
\vspace{-0.4 pc}
\begin{equation}\label{eq:357}
\pi_m:L^1\rightarrow V_m:v\mapsto\sum_{j=1}^{m}{(v,e_j)e_j}.
\vspace{-0.4 pc}
\end{equation}
The previous expression is well defined thanks to the regularity of the family $\{e_i\}_{i\in\mathbb{N}}$; in addition, it permits to apply the projector to functions not belonging to $H$. This will be useful, e.g., in the forthcoming Lemma \eqref{lem:7}.\\
For each $m\in\mathbb{N}$, we look for a stochastic process $\phi_m=\sum_{j=1}^{m}{c_j(t)e_j(x)}$ such that
\vspace{-0.4 pc}
\begin{subequations}
\label{eq:910}
\begin{empheq}[left={}\empheqlbrace]{align}
  & \,\,\m\phi_m=\pi_m(-u\cdot\nabla \phi_m+\Delta\mu_m)\m t+\m w_m, \label{eq:306}\\
  & \,\, \mu_m:=\pi_m(a\phi_m+\phi^3_m-\phi_m-J\ast\phi_m),\nonumber\\
  & \,\,\phi_m(0)=\pi_m\phi_0,\nonumber
\end{empheq}
\end{subequations}
where $w_m:=\pi_m w$ is a $V_m$-valued Wiener process. 
If we take the $H$-inner product of \eqref{eq:306} with $e_1,\cdots,e_m$ we see that the resulting $\mathbb{R}^m$-valued stochastic differential equation has a locally Lipschitz deterministic integrand. Therefore, for each $m\in\mathbb{N}$, Problem \eqref{eq:910} admits a unique solution. This solution is in principle defined up to some random variable $\zeta_m$, hence for all $t\in[0,\zeta_m(\omega))$. See \cite[Theorem 3.1.]{IW} for full details. We will show in the forthcoming \emph{Step 2} that the solutions $\phi_m$ exist (a.s.) for every $t\in[0,T]$.
\vspace{1 pc}\\
\emph{Step 2: Time domain of $\{\phi_m\}_{m\in\mathbb{N}}$ and some preliminary inequalities}.
We now prove that
\begin{equation}\label{eq:343}
\zeta_m(\omega)\geq T,\quad\mbox{ a.s. in }\Omega,\,\,\forall m\in\mathbb{N},
\end{equation} 
which implies that $\phi_m$ is defined on $[0,T]$ (a.s.). We use some ideas from \cite[p.~132, Proof of Theorem 12.1]{RW}, which however cannot be applied directly since condition \cite[p.~132, (12.3)]{RW} is not satisfied for our finite-dimensional stochastic differential equations.

Let us fix $m\in\mathbb{N}$. If we consider the $\mathbb{R}^m$-valued stochastic differential equation associated with $\phi_m$, we see that its determistic integrand $b_m$ and its stochastic integrand $\sigma_m$ 
are locally Lipschitz. In addition, they depend on time exclusively via $\phi_m$. For each $N\in\mathbb{N}$, we can define a couple $\left\{b_{m,N},\sigma_{m,N}\right\}$, agreeing with $\left\{b_m,\sigma_m\right\}$ on $\{\{s,x\}\in\mathbb{R}\times\mathbb{R}^m:s\geq 0,\,\,|x|\leq N\}$, such that $b_{m,N}$ and $\sigma_{m,N}$ are \emph{globally} Lipschitz. As a consequence, \cite[p.~128, Theorem 11.2]{RW} guarantees that there is a unique solution $\phi_{m,N}$ associated to $\left\{b_{m,N},\sigma_{m,N}\right\}$ and defined on $[0,+\infty)$ (a.s.). We define a sequence of stopping times as follows
$$
\tau_N :=\inf\{\tau>0:\|\phi_{m,N}(\tau)\|\geq N\}\wedge N
$$
The sequence $\{\tau_N\}_{N\in\mathbb{N}}$ is obviously increasing. Moreover \cite[p.~131, Corollary 11.10]{RW} implies that 
\begin{equation}\label{eq:345}
\phi_m=\phi_{m,N}\qquad\mbox{on }[0,\tau_N].
\end{equation}
We apply the It\"o formula to the functional $F(\phi_{m,N}(t)):=\|\phi_{m,N}(t)\|^2$ and, for each $t\in[0,\tau_N(\omega)\wedge T)$, we obtain 
the equality 
\begin{eqnarray*}
\|\phi_{m,N}(t)\|^2 & = & \|\phi_{m,N}(0)\|^2+2\int_{0}^{t}{(\phi_{m,N}(s),-u\cdot\nabla\phi_{m,N}(s)+\Delta\mu_{m,N}(s))\m s}\nonumber\\
& + & \int_{0}^{t}{\mbox{tr}(Q_m)}+2\int_{0}^{t}{(\phi_{m,N}(s),\m w_m(s))},
\end{eqnarray*}
where $Q_m$ denotes the covariance operator of $w_m$ and $\mu_{m,N}$ has an obvious meaning. Since $\mbox{div}(u)=0$ and $u=0$ on $\Gamma$ we have $(\phi_{m,N},u\,\cdot\nabla\phi_{m,N})=0$. If we use \cite[p.~436, estimate (4.15)]{CFG} we deduce
\begin{equation}\label{eq:308}
\|\phi_{m,N}(t)\|^2+\int_{0}^{t}{\!\frac{c_0}{2}\|\nabla\phi_{m,N}(s)\|^2\m s}\leq\|\phi_{m,N}(0)\|^2+\int_{0}^{t}{\!\!k\|\phi_{m,N}(s)\|^2\m s}
+ C_1 + 2\int_{0}^{t}{\!\!(\phi_{m,N}(s),\m w_m(s))},
\end{equation}
where $c_0$ is the positive costant from \eqref{eq:353} and $k:=(2/c_0)\|J\|^2_{L^1(\mathbb{R}^d)}$. Therefore, for each $t\in[0,T]$, we have
\begin{eqnarray}\label{eq:319}
& &\mathbb{E}\left[\|\phi_m(t\wedge \tau_N)\|^2\right]+\frac{c_0}{2}\mathbb{E}\left[\int_{0}^{t\wedge\tau_N}{\!\!\!\|\nabla\phi_m(s)\|^2}\right]\leq\mathbb{E}\left[\|\phi_0\|^2\right]+k\mathbb{E}\left[\int_{0}^{t\wedge\tau_N}{\!\!\!\|\phi_m(s)\|^2}\right]+C_2,
\end{eqnarray}
where we have used \eqref{eq:345}. Inequality \eqref{eq:319} implies
\begin{eqnarray*}
\mathbb{E}\left[\|\phi_m(t\wedge \tau_N)\|^2\right]\leq C_2+\mathbb{E}\left[\|\phi_0\|^2\right]+k\!\int_{0}^{t}{\!\mathbb{E}\left[\|\phi_m(s\wedge\tau_N)\|^2\right]}.
\end{eqnarray*}
The Gronwall lemma consequently gives
\begin{equation}\label{eq:346}
\mathbb{E}\left[\|\phi_m(t\wedge\tau_N)\|^2\right]\leq C_3,\qquad\mbox{for }t\in[0,T].
\end{equation}
We notice that the previous inequality holds for every $T>0$. Hence we may take $K=2T$, and use the Markov inequality 
to deduce that, for $N>K$, we have
\begin{eqnarray}\label{eq:344}
\mathbf{m}(\tau_N<K) & \leq &\mathbf{m}(\|\phi_{m,N}(\tau_N\wedge K)\|\geq N)\leq\mathbf{m}(\|\phi_m(\tau_N\wedge K)\|\geq N)\nonumber\\
& \leq & \frac{\mathbb{E}\left[\|\phi_m(K\wedge\tau_N)\|^2\right]}{N^2}\leq\frac{C(K,\phi_0)}{N^2}\rightarrow 0
\end{eqnarray}
for $N\rightarrow+\infty$. Computation \eqref{eq:344} clearly implies that $\mathbf{m}(\sup_{N}\tau_N>T)=1$. Hence \eqref{eq:343} holds. In addition 
\begin{equation}\label{eq:362}
\phi_m\in\mathcal{C}([0,T];\mathcal{C}^{\infty}(\overline{D})),\qquad \forall m\in\mathbb{N}.
\vspace{0.5 pc}
\end{equation}
\emph{Step 3: Main estimates for the family $\{\phi_m\}_{m\in\mathbb{N}}$}. Let $\phi^N_m(t):=\phi_m(t\wedge\tau_N)$. 
Clearly $\phi^{N+1}_m(s)=\phi^{N}_{m}(s)$ for all $0<s<t\wedge\tau_N$ (a.s.). Therefore $t\wedge\tau_N\uparrow t$ (a.s.) and $\phi_m(t)=\lim_{N\rightarrow+\infty}{\phi^N_{m}(t)}$ (a.s.). We use the Fatou Lemma on both the spatial domain and the probability space to deduce that 
\begin{eqnarray}\label{eq:309}
& & \mathbb{E}\left[\|\phi_m\|^2_{L^2([0,T];H)}\right]=\mathbb{E}\left[\int_{0}^{T}{\|\phi_m\|^2}\right]\leq\mathbb{E}\left[\int_{0}^{T}{\liminf_{N}{\|\phi^N_m\|^2}}\right]\nonumber\\
& = & \int_{0}^{T}{\mathbb{E}\left[\liminf_{N}{\|\phi^N_m\|^2}\right]}\leq\int_{0}^{T}{\liminf_{N}{\mathbb{E}\left[\|\phi^N_m\|^2\right]}}\leq C_4,
\end{eqnarray}
where we have used \eqref{eq:346}. We apply the monotone convergence theorem in \eqref{eq:319} and we conclude that
\begin{equation}\label{eq:311}
\mathbb{E}\left[\|\phi_m\|^2_{L^2([0,T];U)}\right]\leq C_5.
\end{equation}
We now take the supremum of \eqref{eq:308} for $0\leq t\leq T$ and get
\begin{equation}\label{eq:347}
 \mathbb{E}\left[\sup_{0\leq t\leq T}{\|\phi_m(t)\|^2}\right]\leq\mathbb{E}\left[\|\phi_m(0)\|^2\right]+C_1
 +  k\mathbb{E}\left[\|\phi_m\|^2_{L^2([0,T];H)}\right]+2\mathbb{E}\left[\sup_{0\leq t\leq T}{\int_{0}^{t}{(\phi_m(s),\m w_m(s))}}\right]\!\!.
\end{equation}
We estimate the last term of the right hand side of \eqref{eq:347} using Doob's submartingale inequality, and get
\begin{eqnarray*}
& & \left\{\mathbb{E}\left[\sup_{0\leq t\leq T}{\int_{0}^{t}{(\phi_m(s),\m w_m(s))}}\right]\right\}^2 \leq \mathbb{E}\left[\sup_{0\leq t\leq T}{\int_{0}^{t}{(\phi_m(s),\m w_m(s))}}\right]^2\\
& \leq & 4\mathbb{E}\left[\left\{\int_{0}^{T}{(\phi_m(s),\m w_m(s))}\right\}^2\right]=4\mathbb{E}\left[\int_{0}^{T}{(Q_m\phi_m,\phi_m)}\right]
 \leq  4\|Q\|\mathbb{E}\left[\int_{0}^{T}{\|\phi_m\|^2}\right]\leq 4C_6\|Q\|,
\end{eqnarray*}
and hence we deduce
\begin{equation}\label{eq:310}
\mathbb{E}\left[\|\phi_m\|^2_{L^\infty([0,T];H)}\right]\leq C_7.
\end{equation}
We now apply the It\"o formula to the functional $Z(\phi_m(t))$, where
$$
Z(\phi):U\rightarrow\mathbb{R}:\phi\mapsto\int_{D}{\left\{a\frac{\phi^2}{2}+\frac{\phi^4}{4}-\frac{\phi^2}{2}\right\}}-\frac{1}{2}(J\ast\phi,\phi).
$$
We defined $Z$ such that its Fr\'echet derivative is $Z_{\phi}(\phi)=a\phi+\phi^3-\phi-J\ast\phi$, hence $\pi_m(Z_{\phi}(\phi))=\mu_m$. Thus, recalling the hypotheses on $Q$, $u$, and the family $\{e_i\}_{i\in\mathbb{N}}$, we obtain that
\begin{eqnarray}\label{eq:1000}
& & Z(\phi_m(t)) = Z(\phi_m(0))+\int_{0}^{t}{(\mu_m,-u\cdot\nabla\phi_m+\Delta\mu_m)}+\int_{0}^{t}{(Z_{\phi}(\phi_m(s)),\m w_m(s))_{U}}\nonumber\\
& + &\frac{1}{2}\left(\int_{0}^{t}{\!\!\!\int_{D}{(3\phi^2_m+a-1)\sum_{i,j=1}^{m}{Q_{ij}e_ie_j}}}-\int_{0}^{t}{\sum_{i=1}^{m}{(J\ast Q_me_i,e_i)}}\right)\nonumber\\
& & \leq Z(\phi_m(0))+\int_{0}^{t}{\left(\|u\|_{L^{\infty}([0,T]\times\Omega)}\|\nabla\mu_m(s)\|\|\phi_m(s)\|-\|\nabla\mu_m(s)\|^2\right)}\nonumber\\
& + & \int_{0}^{t}{(Z_{\phi}(\phi_m(s)),\m w_m(s))_{U}}+C_8\mbox{tr}(Q)\int_{0}^{t}{\left[\|\phi_m(s)\|^2+\|a\|_{L^{\infty}}+\|J\|_{L^1(\mathbb{R}^d)}\right]},
\end{eqnarray}
where $[Q]_{ij}$, $i,j\in\{1,\cdots,m\}$, denote the entries of the matrix representation of $Q_m$ with respect to the basis $\{e_1,\cdots,e_m\}$. Thanks to the Young inequality, computation \eqref{eq:1000} implies
\begin{eqnarray}\label{eq:321}
Z(\phi_m(t))\,& + &\frac{1}{2}\int_{0}^{t}{\|\nabla\mu_m(s)\|^2}\,\,\leq\,\, Z(\phi_m(0))+\int_{0}^{t}{(Z_{\phi}(\phi_m(s)),\m w_m(s))_{U}} + C(T,u,Q)\int_{0}^{t}{\|\phi_m(s)\|^2}\nonumber\\
& + & C_9(K(Q))\int_{0}^{t}{\left[\|\phi_m(s)\|^2+\|a\|_{L^{\infty}}+\|J\|_{L^1(\mathbb{R}^d)}\right]}.
\end{eqnarray}
If we set
$$
\tau^1_N =
     \left\{
     \begin{array}{l}
      \inf\{\tau>0:\|Z_{\phi}(\phi_m(\tau))\|_{U}\geq N\}\quad\mbox{ if }\exists \tau>0:\|Z_{\phi}(\phi_m(\tau))\|_{U}\geq N,\\
     +\infty\quad\mbox{ if }\|Z_{\phi}(\phi_m(\tau))\|_{U}<N,\quad\forall \tau>0,
     \end{array}
     \right.
$$
we may act on \eqref{eq:321} and $\tau^1_N$ similarly to the computations previously done with $\tau_N$. We deduce
\begin{eqnarray}\label{eq:349}
\mathbb{E}\left[Z(\phi_m(t\wedge\tau^1_N))\right]+\frac{1}{2}\mathbb{E}\left[\int_{0}^{t\wedge\tau^1_N}{\|\nabla\mu_m(s)\|^2}\right]\leq\mathbb{E}\left[Z(\phi_m(0))\right]+C_{10}.
\end{eqnarray}
In addition, the regularity of the trajectories of $\phi_m$ highlighted in \eqref{eq:362} implies that $\tau^{1}_N$ eventually coincides with $T$, (a.s.). Computation \cite[p.~1190, (7.6)]{EM} allows to estimate $\mathbb{E}\left[Z(\phi_m(0))\right]$ \emph{uniformly} in $m\in\mathbb{N}$. In addition we can rely on estimates \eqref{eq:311}, \eqref{eq:310} and deduce from \eqref{eq:321} that
\begin{equation}\label{eq:348}
\mathbb{E}\left[\|\phi_m(t\wedge\tau^1_N)\|^4_{L^4}\right]\leq C_{11}.
\end{equation}
We apply again Fatou's Lemma as we have done to deduce \eqref{eq:309} and we obtain
\begin{eqnarray}\label{eq:313}
& & \mathbb{E}\left[\|\phi_m\|^4_{L^4([0,T];L^4)}\right]=\mathbb{E}\left[\int_{0}^{T}{\|\phi_m\|_{L^4}^4}\right]\leq\mathbb{E}\left[\int_{0}^{T}{\liminf_{N}{\|\phi^N_m\|_{L^4}^4}}\right]\nonumber\\
& = & \int_{0}^{T}{\mathbb{E}\left[\liminf_{N}{\|\phi^N_m\|_{L^4}^4}\right]}\leq\int_{0}^{T}{\liminf_{N}{\mathbb{E}\left[\|\phi^N_m\|_{L^4}^4\right]}}\leq C_{12},
\end{eqnarray}
where we have set $\phi^N_m(t):=\phi_m(t\wedge\tau^1_N)$. We have also used \eqref{eq:348}. We can now apply the monotone convergence theorem in \eqref{eq:349} by passing to the limit with respect to $N\rightarrow+\infty$ and deduce
\begin{equation}\label{eq:324}
\mathbb{E}\left[\|\nabla\mu_m\|^2_{L^2([0,T];H)}\right]\leq C_{13}.
\end{equation}
We exploit \cite[p.~1179, Theorem 4.2.]{EM} and we write, for every $v\in V$ and for every $0\leq t_1<t_2\leq T$
\begin{eqnarray}\label{eq:363}
& & \left|\langle \phi_m(t_2)-\phi_m(t_1),v\rangle_{V',V}\right|=\left|(\phi_m(t_2)-\phi_m(t_1),v)\right|\nonumber\\
& =  &\left|\int_{t_1}^{t_2}{(-u\cdot\nabla\phi_m+\Delta\mu_m,v)}+(w_m(t_2)-w_m(t_1),v)\right|\nonumber\\
& \leq & \left(\left[\int_{t_1}^{t_2}{\|u\|_{L^{\infty}([0,T]\times D)}\|\nabla\phi_m\|+\|\nabla\mu_m\|}\right]+|t_2-t_1|^{\frac{2}{5}}\|w_m\|_{\mathcal{C}^{\frac{2}{5}}([0,T];V')}\right)\|v\|_{V}\nonumber\\
& \leq & C_6(u)|t_2-t_1|^{\frac{2}{5}}\cdot\left[\|\phi_m\|_{L^2([0,T];U)}+\|\nabla\mu_m\|_{L^2([0,T];H)}+\|w_m\|_{\mathcal{C}^{\frac{2}{5}}([0,T];V')}\right]\!\|v\|_{V}.
\end{eqnarray}
In addition, \eqref{eq:362} and \eqref{eq:310} imply
\begin{equation}\label{eq:364}
\mathbb{E}\left[\|\phi_m\|^2_{\mathcal{C}([0,T];H)}\right]\leq C_7.
\end{equation}
The combination of \eqref{eq:363} and \eqref{eq:364} allows us to deduce
\begin{equation}\label{eq:314}
\mathbb{E}\left[\|\phi_m\|_{\mathcal{C}^{2/5}([0,T];V')}\right]\leq C_{14}.
\end{equation}
The constants $C_1,\cdots,C_{14}$ are independent of $m$ but may depend on $\phi_0$, $T$, $u$, $J$, $Q$. Inequalities \eqref{eq:311}, \eqref{eq:310}, \eqref{eq:313}, \eqref{eq:314} imply that $\{\mathbf{P}_m\}_{m\in\mathbb{N}}$, the family of the distributions of $\{\phi_m\}_{m\in\mathbb{N}}$ on $\mathscr{Z}$, is \emph{uniformly} concentrated on $\mathscr{U}$.
\vspace{1 pc}\\
\emph{Step 4: Existence of a weak limit.}
Since $\mathscr{U}$ is compactly embedded in $\mathscr{Z}$ as proved in Theorem \eqref{th:20}, we can use a compactness argument relying on the Prohorov Theorem. 
Estimates \eqref{eq:311}, \eqref{eq:310}, \eqref{eq:313}, \eqref{eq:314}, the definition of the $\mathscr{U}$-norm and the $L^p$ embeddings imply that
\begin{equation}\label{eq:354}
\int_{\mathscr{Z}}{\|\phi\|_{\mathscr{U}}\mathbf{P}_m(\m \phi)}\leq C_{15},\quad\forall m\in\mathbb{N},
\end{equation}
where $C_{15}$ is independent of $m$. In addition we have that the sets 
$$
G_n(\mathscr{U}):=\left\{v\in\mathscr{U}\!:\,\,\|v\|_{L^{2}([0,T];U)}\leq n,\,\|v\|_{L^{\infty}([0,T];H)}\leq n,
\,\|v\|_{L^4([0,T];L^4)}\leq n,\|v\|_{\mathcal{C}^{2/5}([0,T];V')}\leq n\right\}
$$
are compact in $\mathscr{Z}$. A full proof of this fact is given in \emph{Step 1} of \emph{Proof of Theorem} \eqref{th:18} below. 
We deduce
\begin{eqnarray*}
& &\mathbf{P}_m(\mathscr{Z}\setminus G_n(\mathscr{U}))=\int_{\mathscr{Z}\setminus G_n(\mathscr{U})}{\mathbf{P}_m(\m \phi)}=\int_{\mathscr{U}\setminus G_n(\mathscr{U})}{\mathbf{P}_m(\m \phi)}\\
& \leq & \frac{1}{n}\int_{\mathscr{U}\setminus G_n(\mathscr{U})}{\|\phi\|_{\mathscr{U}}\mathbf{P}_m(\m \phi)}\leq\frac{1}{n}\int_{\mathscr{Z}}{\|\phi\|_{\mathscr{U}}\mathbf{P}_m(\m \phi)}\leq \frac{C_{15}}{n}\leq\varepsilon,
\end{eqnarray*}
where the last inequality holds if $n$ is large enough. We have verified the validity of the assumptions of the Prohorov theorem. Hence we deduce that there are a (not relabeled) subsequence $\{\mathbf{P}_m\}_{m\in\mathbb{N}}$ and a probability $\mathbf{P}$ defined on $(\mathscr{Z},\mathcal{B}(\mathscr{Z}))$ such that $\{\mathbf{P}_m\}_{m\in\mathbb{N}}$ weakly converges to $\mathbf{P}$.
\vspace{1 pc}\\
\emph{Step 5: Passage to the limit.}
The probability $\mathbf{P}_m$ satisfies 
\begin{eqnarray}\label{eq:322}
\int_{\mathscr{Z}}{\exp\left\{i\langle \phi(0),\xi\rangle_{H^{-\varepsilon},H^{\varepsilon}}+iC(\phi,v)\right\}\mathbf{P}_m(\m\phi)}=\hat{\Xi}(\xi)\hat{\mathbf{W}}_m\!\!\left(-\frac{\partial v}{\partial t}\right)\!,\quad\forall\xi\in V_r,\,\,\forall v\in\mathscr{V}_r,\,\,m\geq r,
\end{eqnarray}
where $\hat{\mathbf{W}}_m$ denotes the characteristic functional of $w_m$. We prove \eqref{eq:322} by \emph{rigorously} replicating the computations we have done in Remark \eqref{rem:16}. Let $r\leq m$ be positive integers. For each $v\in\mathscr{V}_r$, we apply the It\"o formula to the real-valued process $B(\phi_m):=(\phi_m,v)$ and obtain
\begin{eqnarray}\label{eq:340}
0=(\phi_m(0),v(0))+\int_{0}^{T}{\!\!\int_{D}{\pi_m(-u\cdot\nabla\phi_m+\Delta\mu_m)v}}+\int_{0}^{T}{\left(\phi_m,\frac{\partial v}{\partial t}\right)}-\int_{0}^{T}{\left(\frac{\partial v}{\partial t},w_m\right)}.
\end{eqnarray}
Since $v$ takes values in $V_m$, we can integrate by parts and rewrite \eqref{eq:340} as
\begin{eqnarray*}
0=(\phi_0,v(0))+\int_{0}^{T}{\!\!\int_{D}{\left\{-uv\cdot\nabla\phi_m+\mu\Delta v\right\}}}+\int_{0}^{T}{\left(\phi_m,\frac{\partial v}{\partial t}\right)}-\int_{0}^{T}{\left(\frac{\partial v}{\partial t},w_m\right)}.
\end{eqnarray*}
Hence we can rearrange the terms and add the term $\langle\phi_m(0),\xi\rangle_{H^{-\varepsilon},H^{\varepsilon}}$ in the last equality to obtain
\begin{eqnarray}\label{eq:342}
\langle\phi_m(0),\xi\rangle_{H^{-\varepsilon},H^{\varepsilon}}+C(\phi_m,v)=(\phi_0,\xi)-\int_{0}^{T}{\!\!\left(w_m,\frac{\partial v}{\partial t}\right)},
\end{eqnarray}
where we used the injections $H^{\varepsilon}\hookrightarrow H\hookrightarrow H^{-\varepsilon}$ and the fact that $\xi\in V_m$. Moreover, $w_m$ and $\phi_m(0)$ are independent random variables because of the independence of $w$ and $\phi_0$. We multiply \eqref{eq:342} by $i$, apply the exponential function, take the expected value, use the independence of $w_m$ and $\phi_m(0)$, and \eqref{eq:322} follows. \\
In Remark \eqref{rem:10} we have seen that $\phi\mapsto C(\phi,v)$ is continuous on $\mathscr{Z}$ for every $v\in\mathscr{V}$. In addition, we recall that $w_m\rightarrow w$ in $L^2((\Omega,\mathcal{F},\mathbf{m});\mathcal{C}([0,T];H))$. See \cite[p.~13, Proposition 2.1.10]{PR}. Therefore we can rely on the weak convergence of the sequence $\{\mathbf{P}_m\}_{m\in\mathbb{N}}$, hence take the limit for $m\rightarrow+\infty$ in \eqref{eq:322} and deduce that \eqref{eq:323} holds for $\mathbf{P}$ and $\xi\in\cup_{r=1}^{\infty}{V_r}$, $v\in\cup_{r=1}^{\infty}{\mathscr{V}_r}$. By the density of $\cup_{r=1}^{\infty}{V_r}$ in $H^{\varepsilon}$ (see Remark \eqref{rem:21}) and the density of $\cup_{r=1}^{\infty}{\mathscr{V}_r}$ in $\mathscr{V}$, and using Lebesgue's dominated convergence theorem, 
we deduce that $\mathbf{P}$ satisfies \eqref{eq:323} for each $\xi\in H^{\varepsilon}$ and $v\in\mathscr{V}$. To prove that $\mathbf{P}$ is a weak statistical solution to Problem \eqref{eq:303} in the sense of Definition \eqref{def:2}, it remains to show that $\mathbf{P}$ is concentrated on $\mathcal{B}_{\mathscr{Z}}(\mathscr{U})$.
\vspace{1 pc}\\
\emph{Step 6: $\mathbf{P}$ is concentrated on $\mathcal{B}_{\mathscr{Z}}(\mathscr{U})$.}
It is straightforward to notice that
\begin{equation}\label{eq:350}
\mathscr{U}=\bigcup_{n=1}^{\infty}{G_n(\mathscr{U})}.
\end{equation}
It follows from \eqref{eq:350} that $\mathscr{U}\in\mathcal{B}_{\mathscr{Z}}(\mathscr{U})$. Moreover, as said in \emph{Step 4}, the sets $G_n(\mathscr{U})$'s are compact (hence closed) in $\mathscr{Z}$. As a consequence of \emph{Step 4}, for each $\varepsilon>0$ there is $n\in\mathbb{N}$ such that $\mathbf{P}_m(G_n(\mathscr{U}))\geq 1-\varepsilon$ for each $m\in\mathbb{N}$. We use the Portmanteau Theorem to deduce that $\mathbf{P}(G_n(\mathscr{U}))\geq 1-\varepsilon$. Hence $\mathbf{P}(\mathscr{U})=1$. 
\end{Proof}
\begin{rem}
The definition of the test function space $\mathscr{V}$ is required to use the It\"o formula in \emph{Step 5} of the previous proof. In particular, the requirement on the time double derivative is crucial. 
\end{rem}
\begin{rem} Estimate \eqref{eq:324} gives an estimate on some $H$-projections $\mu_m$ of the chemical potential $\mu$. In Sections \ref{sec:41} and \ref{sec:42}, however, we will need to act on the ``full'' chemical potential $\mu$. For this purpose, we need the following lemma.
\end{rem}
\begin{lemma}\label{lem:7}
With the notation of Theorem \eqref{th:19}, we have
\begin{equation}\label{eq:335}
\int_{\mathscr{Z}}{\|\nabla\mu(\phi)\|^{2}_{L^2([0,T];H)}\mathbf{P}(\mbox{\emph{d}}\phi)}\leq C_{13},
\end{equation}
\vspace{-1 pc}
\begin{equation}\label{eq:336}
\int_{\mathscr{Z}}{\|\phi\|^{2}_{L^2([0,T];U)}\mathbf{P}(\mbox{\emph{d}}\phi)}\leq C_{5},
\end{equation}
where $\mu(\phi):=a\phi-J\ast\phi+\phi^3-\phi$. The integrands are understood to assume the value $+\infty$ whenever $\nabla\mu(\phi)\notin L^2([0,T];H)$ or $\phi\notin L^2([0,T];U)$. 
\end{lemma}
\begin{Proof}
Estimate \eqref{eq:324} implies
\begin{equation}\label{eq:332}
\int_{\mathscr{Z}}{\|\nabla\mu_r(\phi)\|^2_{L^2([0,T];H)}\mathbf{P}_m(\m \phi)}\leq C_{13},\quad \mbox{ for } r\leq m.
\end{equation}
For each $r\in\mathbb{N}$, we define
$$
\Phi_r:\mathscr{Z}\rightarrow\mathbb{R}\cup\{+\infty\}:\phi\mapsto\|\nabla\mu_r(\phi)\|^2_{L^2([0,T];H)}.
$$
We use the definition of $\pi_r$ given in \emph{Proof of Theorem} \eqref{th:19}, \emph{Step 1}, since $\mu(\phi)$ takes values in $L^{4/3}$, a.e. $t\in[0,T]$, for $\phi\in\mathscr{Z}$. The functional $\Phi_r$ is lower semicontinuous in $\mathscr{Z}$. To prove this, let $\phi_n\rightarrow \phi$ in $\mathscr{Z}$. Then, for any (not relabeled) subsequence, there is another (not relabeled) subsequence $\phi_n$ such that $\phi_n(t)\rightarrow \phi(t)$ in $L^4$ for a.e. $t\in[0,T]$. By taking the $H$-inner product of $\phi_n(t)-\phi(t)$ with $e_1,\cdots\!,e_m$ and using the H\"older inequality along with the regularity of $e_1,\cdots\!,e_m$, we deduce that $\mu_r(\phi_n)(t)\rightarrow\mu_r(\phi)(t)$ in $H$. Hence, relying on the equivalence of norms of finite-dimensional vector spaces, we deduce that $\nabla\mu_r(\phi_n)(t)\rightarrow \nabla\mu_r(\phi)(t)$ in $H$ for a.e. $t\in[0,T]$. Since the original subsequence $\phi_n$ was arbitrary, we may apply the Fatou Lemma to deduce that $\Phi_r$ is lower semicontinuous in $\mathscr{Z}$. In addition $\Phi_r$ is nonnegative and thus is trivially bounded from below. Hence we may apply the Portmanteau Theorem and write
\begin{equation}\label{eq:333}
\int_{\mathscr{Z}}{\Phi_r(\phi)\mathbf{P}(\m \phi)}\leq\liminf_{m}{\int_{\mathscr{Z}}{\Phi_r(\phi)\mathbf{P}_m(\m \phi)}}\leq C_{13},
\end{equation}
where the last inequality follows from \eqref{eq:332}. But now, thanks to the monotonicity of the $H$-norm under projection to growing subspaces, we may use the monotone convergence Theorem in \eqref{eq:333} to deduce
\eqref{eq:335}. The proof of \eqref{eq:336} is similar. It suffices to consider $\Phi_r(\phi)=\|\pi_r \phi\|_{L^2([0,T];U)}$. A finite-dimensional argument identical to the one used before shows that $\Phi_r$ is lower semicontinuous in $\mathscr{Z}$ and bounded from below. Hence we once again apply the Portmanteau theorem and get
\begin{equation}\label{eq:355}
\int_{\mathscr{Z}}{\Phi_r(\phi)\mathbf{P}(\m \phi)}\leq\liminf_{m}{\int_{\mathscr{Z}}{\Phi_r(\phi)\mathbf{P}_m(\m \phi)}}\leq C_{5},
\end{equation}
where we have used \eqref{eq:311}. Since $\{e_i\}_{i\in\mathbb{N}}$ is an orthogonal basis in $U$ as well, we may use a similar monotonicity argument with respect to the $U$-norm under projection to growing subspaces, and we may replicate the application of the monotone convergence Theorem in \eqref{eq:355} to deduce \eqref{eq:336}.
\end{Proof}
\begin{rem}\label{rem:22}
In the previous lemma, we have implicitly used the fact that, for every $x\in L^{4/3}$: \\
$\mbox{ }\,\,\,$- the sequence $\{\|\nabla(\pi_rx)\|\}_{r\in\mathbb{N}}$ converges to $\|\nabla x\|$ if $x\in U$, and diverges to $+\infty$ if $\nabla x\notin H$.\\
$\mbox{ }\,\,\,$- the sequence $\{\|\pi_rx\|_{U}\}_{r\in\mathbb{N}}$ converges to $\|x\|_{U}$ if $x\in U$, and diverges to $+\infty$ if $x\notin U$.
\end{rem}
\begin{rem}\label{rem:14} Let $\mathbf{P}$ be the weak statistical solution built in Theorem \eqref{th:19}. Lemma \eqref{lem:7} and Remark \eqref{rem:22} imply that 
$$
\mathbf{P}(\{v\in\mathscr{Z}:\nabla\mu(v)\in L^2([0,T];H)\})=1,
\vspace{-0.4 pc}
$$
$$
\mathbf{P}(\{v\in\mathscr{Z}:\nabla (av-v)\in L^2([0,T];H)\})=1,
\vspace{-0.4 pc}
$$ 
$$
\mathbf{P}(\{v\in\mathscr{Z}:\nabla(J\ast v)\in L^2([0,T];H)\})=1.
$$ 
Hence, recalling \emph{Proof of Theorem} \eqref{th:19}, \emph{Step 6}, we deduce that the set
$$
\bigl\{v\in\mathscr{U}:\,\,\nabla(av-v)\in L^2([0,T];H),\nabla(J\ast v)\in L^2([0,T];H),\,\,3v^2\nabla v\in L^2([0,T];H)\bigr\}
$$
has $\mathbf{P}$-probability one. This fact will be useful in Section \ref{sec:42}.
\end{rem}

\section{Strong solutions}\label{sec:41}

The main purpose of this brief section is to outline the strategy of the proof of the existence and uniqueness of \emph{strong solutions} to Problem \eqref{eq:303}. To be able to do this we must first explain what we mean by strong solution.

\begin{mydef}\label{def:4}
A process $\phi=\phi(t,x,\omega)$ defined on the probability space $(\Omega,\mathcal{F},\mathbf{m})$ is a {\bfseries strong} solution to Problem \eqref{eq:303} if
\begin{enumerate}
\item $\phi$ satisfies 
\begin{equation}\label{eq:331}
\mathscr{D}(\phi(\omega))=\left\{ \phi_0(\omega),\frac{\partial w}{\partial t}(\omega)\right\}\qquad \mathbf{m}-\mbox{a.s.}
\end{equation}
\item the mapping $\omega\mapsto \phi(\omega)$ is a random variable from $(\Omega,\mathcal{F})$ to $(\mathscr{U}_1,\mathcal{B}_{\mathscr{Z}}(\mathscr{U}_1))$, where we have set
$$
\mathscr{U}_1:=\{v\in\mathscr{U}:\,\,\nabla\mu(v)\in L^2([0,T];H)\}.
$$ 
\end{enumerate}
\end{mydef}
The above definition deserves some comments.\vspace{1 pc}\\
{\bfseries Comment 1}.
Our definition of strong solution to Problem \eqref{eq:303} is weaker than the definition of \emph{classical strong solution} for a stochastic partial differential equation (cf. for instance \cite[p.~73, Definition 4.2.1.]{PR}). As a matter of fact, our solution is a \emph{variational} solution, hence it is paired to the space of \emph{test functions} $\mathscr{V}$. In addition, even though there is no restriction on the final positive time $T$ (we did not have to make any assumptions up to now and the same will happen in the sequel), our solution satisfies a relation which involves the \emph{entire} time interval $[0,T]$ instead of any arbitrary interval $[0,t]$, for $t\in[0,T]$.\\
The motivation for giving the above definition of strong solution is the relative lack of regularity of the equation (lack of sublinear growth conditions in particular) which does not allow to apply the classical theorems of existence and uniqueness (cf. for instance \cite[p.~75, Theorem 4.2.4.]{PR}).\vspace{1 pc}\\
{\bfseries Comment 2}. A strong solution $\phi$ is required to take values in a regular subset of $\mathscr{U}$, namely in $\mathscr{U}_1$. The reason for requiring this is suggested, in particular, by estimate \eqref{eq:324}, Lemma \eqref{lem:7} and Remark \eqref{rem:14}. As a matter of fact, it stands to reason that the \emph{distribution} of the weak solution $\mathbf{P}$ built in Theorem \eqref{th:19} might give indications on where one can look for ``stronger'' solutions. Since $\mathbf{P}$ gives full probability to the set 
$$
\{v\in\mathscr{U}:\nabla\mu(v)\in L^2([0,T];H)\},
$$
we consequentially include the above set in the definition.\vspace{1 pc}\\
{\bfseries Comment 3}. We have introduced the operator $\mathscr{D}$ in Remark \eqref{rem:13}, item (b). This very operator encodes the initial condition $\phi_0$ and the ``deterministic body'' of our Cahn-Hilliard equation. In equation \eqref{eq:331} the solution $\phi$ is mapped via $\mathscr{D}$ to the initial condition $\phi_0$ and to the given white noise $\partial w/\partial t$. The functional equality \eqref{eq:331}, which is justified by \eqref{eq:338}, is our rewriting of Problem \eqref{eq:303}.\\

We may now outline the strategy we will use to prove the existence and uniqueness of strong solutions, keeping in mind Comments 2, 3, and with the help of the results proved in Section \ref{sec:40}.
\begin{itemize}
\item We first prove that the operator $\mathscr{D}$, if restricted to $\mathscr{U}_1$, is injective. This fact follows from a uniqueness result which will also be proved. 
\item We then prove that the above restricted operator has the initial condition $\phi_0$ and the white noise $\partial w/\partial t$ in its image with full $\mathbf{m}$-probability. This observation is \emph{crucial}, because it allows us to perform a simple inversion in order to construct a random variable $\phi$ (a strong solution) satisfying \eqref{eq:331}. Because the restricted operator is injective and satisfies specific measurability properties, existence and uniqueness of a strong solution are deduced. For this argument, we use in particular Theorem \eqref{th:19}, Lemma \eqref{lem:70} and Remark \eqref{rem:14}.
\end{itemize}


\section{Existence and uniqueness of a strong solution}\label{sec:42}

In this final section we provide the proofs which have been outlined at the end of the previous section. We first state and prove a result which will be used to deduce the uniqueness of a strong solution.
\begin{theorem}\label{th:17}
Let $\phi_1$, $\phi_2$ be two strong solutions to Problem \eqref{eq:303} (for the same $\phi_0\in U$, i.e. $\phi_1(0)=\phi_2(0)=\phi_0$) in the sense of Definition \eqref{def:4}. Then
$$
\phi_1(t)=\phi_2(t)\quad\mbox{in }U',\quad\mbox{for a.e. }t\in[0,T].
$$
\end{theorem}
\begin{rem}\label{rem:17}
The Proof of Theorem \eqref{th:17} is carried out by means of purely deterministic arguments. The stochastic noise of our version of the Cahn-Hilliard equation is additive and its stochastic integrand is constant (the identity operator); hence, when we subtract the expressions associated with two strong solutions, the stochastic noise vanishes from the computations. If the stochastic integrand were not constant, the proof of the uniqueness would be significantly more complicated and we would be forced to rely on a radically different theory: in fact, the stochasticity could not be removed. 
\end{rem}
\begin{Proof}[Proof of Theorem \eqref{th:17}]
Let $r:=\phi_1-\phi_2$. In analogy with \cite[p.~1195, computations (8.1)-(8.3)]{EM}, we can define, for every $\xi\in V_j$ and $t\in[0,T]$, the sequence 
$$
v_n(s):=
\left\{
 \begin{array}{l}
     \displaystyle (t-s)\xi,\quad\mbox{if }s\in\left[0,t-1/(4n^2)\right]\!,\vspace{0.3 pc}\\
     \left(n(s-[t+1/(4n^2)])\right)^2\xi,\quad\mbox{if }s\in\left[t-1/(4n^2),\,t+1/(4n^2)\right]\!,\vspace{0.3 pc}\\
     0,\quad\mbox{if }s\in(t+1/(4n^2),T].
     \end{array}
     \right.
$$
Here $j$ is an arbitrary nonnegative integer. As a result, $v_n\in \mathscr{V}$ and
$$
\frac{\partial v_n}{\partial s}=
\left\{
 \begin{array}{l}
     -\xi,\quad\mbox{if }s\in\left[0,t-1/(4n^2)\right]\!,\vspace{0.3 pc}\\
     2n^2(s-[t+1/(4n^2)])\xi,\quad\mbox{if }s\in\left[t-1/(4n^2),\,t+1/(4n^2)\right]\!,\vspace{0.3 pc}\\
     0,\quad\mbox{if }s\in(t+1/(4n^2),T].
     \end{array}
     \right.
$$ We now evaluate $C(\phi_1,v_n)-C(\phi_2,v_n)$. We can write
\begin{eqnarray}
& 0 &=-\int_{0}^{T}{\!\!\int_{D}{ru\cdot\nabla v_n}}-\int_{0}^{T}{\!\!\int_{D}{(\mu_1-\mu_2)}\Delta v_n}-\int_{0}^{T}{\!\!\left(r,\frac{\partial v_n}{\partial s}\right)}\nonumber\\
& = & -\int_{0}^{t-1/(4n^2)}{\!\!\int_{D}{ru\cdot(t-s)\nabla \xi}}-\int_{0}^{t-1/(4n^2)}{\!\!\int_{D}{(\mu_{1}-\mu_{2})}(t-s)\Delta \xi}\nonumber\\
& - & \int_{t-1/(4n^2)}^{t+1/(4n^2)}{\!\!\int_{D}{ru\cdot(n(s-[t-1/(4n^2)]))^2\nabla \xi}}\nonumber\\
& - & \int_{t-1/(4n^2)}^{t+1/(4n^2)}{\!\!\int_{D}{(\mu_{1}-\mu_{2})}(n(s-[t-1/(4n^2)]))^2\Delta \xi}\nonumber\\
& + & \int_{0}^{t-1/(4n^2)}{\!\!\left(r,\xi\right)}-\int_{t-1/(4n^2)}^{t+1/(4n^2)}{(r,2n^2(s-[t+1/(4n^2)])\xi)},\label{eq:367}
\end{eqnarray}
where $\mu_1:=\mu(\phi_1)$, $\mu_2:=\mu(\phi_2)$. We use the Lebesgue dominated convergence theorem in \eqref{eq:367} and deduce
\begin{eqnarray}\label{eq:358}
-\int_{0}^{t}{\!\!\int_{D}{ru\cdot(t-s)\nabla \xi}}+\int_{0}^{t}{\!\!\int_{D}{\nabla(\mu_{1}-\mu_{2})}\cdot(t-s)\nabla\xi}+\int_{0}^{t}{\!\!\left(r,\xi\right)}=0.
\end{eqnarray}
We have also used integration by parts to treat the term $(\mu_1-\mu_2,\Delta v)$. If we differentiate \eqref{eq:358} with respect to $t$, we obtain
\begin{eqnarray*}
\langle r(t),\xi\rangle_{V',V}-\int_{0}^{t}{\!\!\int_{D}{ru\cdot\nabla\xi}}+\int_{0}^{t}{\!\!\int_{D}{(\nabla\mu_{1}-\nabla\mu_{2})\cdot\nabla\xi}}=0,\qquad\forall\xi\in V_j.
\end{eqnarray*}
Since $j$ is an arbitrary nonnegative integer, we may use the density of $\cup_{j=1}^{\infty}V_j$ in $V$ and deduce that
$$
\langle r(t),\xi\rangle_{V',V}-\int_{0}^{t}{\!\!\int_{D}{ru\cdot\nabla\xi}}+\int_{0}^{t}{\!\!\int_{D}{(\nabla\mu_1-\nabla\mu_2)\cdot\nabla\xi}}=0,\qquad\forall\xi\in V.
$$
Thanks to the density of $V$ in $U$, and the definition of $\mathscr{U}_1$ we get
\begin{equation}\label{eq:315}
\langle r(t),\xi\rangle_{U',U}-\int_{0}^{t}{\!\!\int_{D}{ru\cdot\nabla\xi}}+\int_{0}^{t}{\!\!\int_{D}{(\nabla\mu_1-\nabla\mu_2)\cdot\nabla\xi}}=0,\quad\forall\xi\in U.
\end{equation}
Equation \eqref{eq:315} implies that $\partial r/\partial t\in L^2([0,T];U')$ since, for every $\xi\in U$ and $\phi\in\mathcal{C}^{\infty}_{0}[0,T]$
\begin{eqnarray*}
& & \int_{0}^{T}{\langle r(t),\xi\rangle_{U',U}\phi'(t)}=\int_{0}^{T}{\!\!\left[\int_{0}^{t}{(ru,\nabla\xi)}-\int_{0}^{t}{(\nabla\mu_1-\nabla\mu_2,\nabla\xi)}\right]\phi'(t)}\\
& = & -\int_{0}^{T}{\left[(ru,\nabla\xi)-(\nabla\mu_1-\nabla\mu_2,\nabla\xi)\right]\phi(t)}=-\int_{0}^{T}{\left\langle \frac{\partial r}{\partial t},\xi\right\rangle_{U',U}\phi(t)},
\end{eqnarray*}
where we have defined $\langle \partial r/\partial t,\xi\rangle_{U',U}:=(ru-(\nabla\mu_1-\nabla\mu_2),\nabla\xi)$. Hence equation \eqref{eq:315} leads to
\begin{equation}\label{eq:326}
\left\langle \frac{\partial r}{\partial t}, \xi\right\rangle_{U',U}+(\nabla\mu_1-\nabla\mu_2,\nabla\xi)=(ru,\nabla\xi),\quad \xi\in U.
\end{equation}
We can therefore act as in \cite[Section 4, Proposition 5]{FG}, from which we reproduce only the computations we need. Since $r(0)=\phi_1(0)-\phi_2(0)=0$, it is clear that $(r(t),1)=0$. We consider the operator $B:\mathcal{D}(B)=V\rightarrow \tilde{H}:u\mapsto-\Delta u,$ where $\tilde{H}:=\{u\in H:(u,1)=0\}$. Then, if we take $\xi=B^{-1}r(t)\in \mathcal{D}(B)$, for almost every $t\in[0,T]$, \eqref{eq:326} implies
$$
\frac{\m}{\m t}\| B^{-1/2}r\|^2+2(\mu_1-\mu_2,r)=2(u,r\nabla B^{-1}r).
$$
If we apply Lagrange's theorem to $F''$ and use \eqref{eq:353}, we obtain
\begin{equation}\label{eq:327}
\frac{\m}{\m t}\| B^{-1/2}r\|^2+2c_0\|r\|^2\leq 2(J\ast r,r)+C\|u\|_{L^{\infty}([0,T]\times D)}\|r\|\| B^{-1/2}r\|.
\end{equation}
Moreover, assumption (iii) and the Young inequality imply
\begin{equation}\label{eq:328}
|(J\ast r,r)|\leq\| B^{1/2}(J\ast r)\|\| B^{-1/2}r\|\leq\frac{c_0}{4}\|r\|^2+C'\| B^{-1/2}r\|^2.
\end{equation}
where $C'>0$ depends on $c_0$ and $J$. If we combine \eqref{eq:327}-\eqref{eq:328} and use once again the Young inequality to control the last term of the right hand side of \eqref{eq:327} we obtain that $r\equiv 0$.
\end{Proof}
\begin{rem}
In the \emph{Proof of Theorem} \eqref{th:17} we have used only certain assumptions contained in the statement of \cite[Section 4, Proposition 5]{FG} because we only needed to prove uniqueness. In particular, we have not made any additional requirement upon the spatial dimension $d$.
\end{rem}
As an immediate consequence we get the following corollary.
\begin{corol}\label{cor:2}
The (restricted) operator $\mathscr{D}|_{\mathscr{U}_1}:\mathscr{U}_1\rightarrow H^{-\varepsilon}\times\mathscr{V}':\phi\mapsto \{\phi(0),C(\phi)\}$ 
is injective.
\end{corol}
We now state and prove our final theorem. 
\begin{theorem}\label{th:18}
Let $d\leq 3$. Let $w$ be a $H$-valued $Q$-Wiener process and let $w$, $u$, $Q$, $J$, $F$, $\phi_0$ and $\{e_j\}_{j\in\mathbb{N}}$ satisfy properties \emph{(i)-(viii)} listed in Section \ref{sec:5}. Let $\phi_0$ be a $U$-valued random variable such that 
$$
\mathbb{E}\left[\|\phi_0\|_{U}^2+\int_{D}{\frac{\phi^4_0}{4}}-\int_{D}{\frac{\phi^2_0}{2}}\right]<+\infty.
$$ 
Then Problem \eqref{eq:303} admits a unique strong solution (in the sense that two strong solutions coincide for all $\omega\in\Omega$ except for a set of $\mathbf{m}$-measure zero).
\end{theorem}
\begin{Proof}
The hypotheses of Theorem \eqref{th:19} being satisfied, we have the weak solution $\mathbf{P}$ 
constructed in \emph{Proof of Theorem} \eqref{th:19}.
\vspace{1 pc}\\
\emph{Step 1: Construction of suitable $\mathscr{Z}$-compact sets.}
Let us consider the countable family of the sets
\begin{eqnarray*}
& & C_{j} := \bigl\{v\in\mathscr{Z}:\|v\|_{L^2([0,T];U)}\leq j,\,\,\|v\|_{L^{\infty}([0,T];H)}\leq j,\,\,\|v\|_{L^4([0,T];L^4)}\leq j,\\
& &\|v\|_{\mathcal{C}^{2/5}([0,T];V')}\leq j,\,\,\|\nabla(J\ast v)\|_{L^2([0,T];H)}\leq j,\,\,\|3v^2\nabla v\|_{L^2([0,T];H)}\leq j,\|\nabla(av-v)\|_{L^2([0,T];H)}\leq j\bigr\},
\end{eqnarray*}
indexed by $j\in\mathbb{N}$. We show that $C_j$ is a compact set in $\mathscr{Z}$. Let $v_n$ be an arbitrary sequence in $C_j$. Because of the compact injection $\mathscr{U}\hookrightarrow\mathscr{Z}$, there is a (not relabeled) subsequence $v_n\rightarrow v$ in $\mathscr{Z}$. We show that $v\in C_j$. Since $v_n\rightarrow v$ in $\mathcal{C}([0,T];H^{-\varepsilon})$, we deduce that $v_n\rightarrow v$ in $\mathcal{C}([0,T];V')$. It is also obvious that
$$
\|v_n\|_{\mathcal{C}([0,T];V')}+\frac{\|v_n(t_1)-v_n(t_2)\|_{V'}}{|t_1-t_2|^{2/5}}\leq j,\quad\forall n\in\mathbb{N},\quad\forall t_1,t_2:\,0\leq t_1<t_2\leq T.
$$ 
If, in the previous inequality, we first take the limit for $n\rightarrow +\infty$ and then the supremum for all $0\leq t_1<t_2\leq T$ we deduce
$$
\|v\|_{\mathcal{C}^{2/5}([0,T];V')}\leq j.
$$ 
Since $v_n$ is bounded in $L^2([0,T];U)$ and $L^4([0,T];L^4)$, we deduce that $v$ is the weak limit both in $L^2([0,T];U)$ and $L^4([0,T];L^4)$ for a (not relabeled) subsequence $v_n$. Hence
\vspace{-0.2 pc}
$$
\|v\|_{L^2([0,T];U)}\leq j,
\vspace{-0.2 pc}
$$
$$
\|v\|_{L^4([0,T];L^4)}\leq j.
$$
We may use a similar argument for the weak convergences in $L^p([0,T];H)$ and deduce
$$
\|v\|_{L^{\infty}([0,T];H)} = \lim_{p\rightarrow+\infty}\|v\|_{L^p([0,T];H)}\leq\lim_{p\rightarrow+\infty}\{\liminf_{n}{\|v_n\|_{L^p([0,T];H)}}\}\leq\lim_{p\rightarrow+\infty}T^{1/p}j=j.
$$
It follows that $v\in G_{j}(\mathscr{U})$. Moreover, the computations done so far in this proof show that $G_{j}(\mathscr{U})$ is a compact set in $\mathscr{Z}$. This fact was used in \emph{Proof of Theorem} \eqref{th:19}, \emph{Steps 4 and 6}.\\
We rely on the $[0,T]$-almost sure convergence in $L^4$ previously noticed to apply Fatou's lemma and deduce
$$
\|\nabla(J\ast v)\|_{L^2([0,T];H)}=\|\nabla(J)\ast v\|_{L^2([0,T];H)}\leq j.
$$
Since $v_n\rightharpoonup v$ in $L^2([0,T];U)$, we deduce that
$$
\|\nabla(av-v)\|_{L^2([0,T];H)}\leq j,
$$
thanks to the regularity of $a$. We now turn to the last and most delicate term appearing in the definition of $C_j$. Let us consider a (not relabeled) subsequence $v_n\rightharpoonup v$ in $L^2([0,T];U)$ such that $v_n(t)\rightarrow v(t)$ in $L^4$ for a.e. $t\in[0,T]$. For each $k\in\mathbb{N}$ and $n\in\mathbb{N}$, we define some suitable \emph{truncated} functions as follows
$$
z^k_n:=\min\{3v^2_n;k\},
$$
\vspace{-1 pc}
$$
z^k:=\min\{3v^2;k\}.
$$
The sequence $\{z^k_n\}_{n\in\mathbb{N}}$ is clearly bounded in $L^{\infty}([0,T];L^{\infty})$. 
Using the time almost sure convergence $v_n(t)\rightarrow v(t)$ in $L^4$ and the H\"older inequality we deduce that, for a.e. $t\in[0,T]$
\begin{eqnarray}\label{eq:359}
& &\|z^k_n-z^k\|^2=\int_{D}{\left|\min\{3v^2_n;k\}-\min\{3v^2;k\}\right|^2}\leq9\int_{D}{|v^2_n-v^2|^2}\nonumber\\
& = & 9\int_{D}{|v_n-v|^2|v_n+v|^2} \leq 36\|v_n-v\|^2_{L^4}\!\left\{\|v_n\|^2_{L^4}+\|v\|^2_{L^4}\right\}\rightarrow 0
\end{eqnarray}
as $n\rightarrow+\infty$. In computations \eqref{eq:359} we have also used to elementary inequalities
$$
\left|\min\{a;c\}-\min\{b;c\}\right|\leq |a-b|,\quad\forall a,b,c\in[0,+\infty),
$$
\vspace{-1.3 pc}
$$
(a+b)^p\leq a^p+b^p,\quad\forall a,b\in[0,+\infty),\,p\in(0,1).
$$
\vspace{-0.8 pc}\\
In addition, $\|z^k_n-z^k\|^2\leq |D|k^2$ for every $t\in[0,T]$ and every $n\in\mathbb{N}$. Hence Lebesgue's dominated convergence theorem implies that $z^k_n\rightarrow z^k$ in $L^2([0,T];H)$ as $n\rightarrow +\infty$. 
Moreover, $\{z^k_n\nabla v_n\}_{n\in\mathbb{N}}$ is bounded in $L^2([0,T];H)$, hence a (not relabeled) subsequence $z^k_n\nabla v_n\rightharpoonup h$ in $L^2([0,T];H)$ as $n\rightarrow+\infty$.
In addition $\nabla v_n\rightharpoonup \nabla v$ in $L^2([0,T];H)$ since $v_n\rightharpoonup v$ in $L^2([0,T];U)$. Hence, for each $\ell\in L^{\infty}([0,T];L^{\infty})$, we deduce
\begin{eqnarray*}
& & \left|(z^k_n\nabla v_n-z^k\nabla v\,,\,\ell)_{L^2([0,T];H)}\right|\leq\left|\int_{0}^{T}{(z^k_n-z^k)\nabla v_n\cdot\ell}\right|
+\left|\int_{0}^{T}{z^k(\nabla v_n-\nabla v)\cdot\ell}\right|\\
&\leq&\|\ell\|_{L^{\infty}([0,T];L^{\infty})}\|z^k_n-z^k\|_{L^2([0,T];H)}\|\nabla v_n\|_{L^2([0,T];H)}+\left|\int_{0}^{T}{z^k(\nabla v_n-\nabla v)\cdot\ell}\right|\rightarrow 0,
\end{eqnarray*}
as $n\rightarrow +\infty$. Because of the density of $L^{\infty}([0,T];L^{\infty})$ in $L^2([0,T];H)$, we deduce that $h= z^k\nabla v$. We can hence rely on the lower semicontinuity property for weakly convergent sequences and deduce
\begin{equation}\label{eq:334}
\int_{0}^{T}{\|z^k\nabla v\|^2}\leq\left[\liminf_{n}\left(\int_{0}^{T}{\|z^k_n\nabla v_n\|^2}\right)^{\!1/2}\right]^{2}\leq\left[\liminf_{n}\left(\int_{0}^{T}{\|3v^2_n\nabla v_n\|^2}\right)^{\!1/2}\right]^{2}\leq j^2.
\end{equation}
For a.e. $t\in[0,T]$ we have that $z^k\rightarrow 3v^2$ a.e. in $D$ for $k\rightarrow +\infty$. Hence, applying Fatou's lemma
in space and time and using \eqref{eq:334}, we deduce
$$
\int_{0}^{T}{\|3v^2\nabla v\|^2}\leq\int_{0}^{T}{\liminf_{k}{\|z^k\nabla v\|^2}}\leq\liminf_{k}\int_{0}^{T}{\|z^k\nabla v\|^2}\leq j^2
$$
and consequently 
$$
\|3v^2\nabla v\|_{L^2([0,T];H)}=\left(\int_{0}^{T}{\|3v^2\nabla v\|^2}\right)^{1/2}\leq j.
$$
We conclude that $C_j$ is compact in $\mathscr{Z}$.
\vspace{1 pc}\\
\emph{Step 2: Costruction of a suitable restriction of $\mathscr{D}$}.
Set $X:=\cup_{j\in\mathbb{N}}{C_{j}}\in\mathcal{B}_{\mathscr{Z}}(\mathscr{U})$. Thanks to Remark \eqref{rem:14} we deduce that $\mathbf{P}(X)=1$. \\
Since any $v\in C_j$ takes values in $U$, it is easy to verify that $\nabla(v^3)=3v^2\nabla v$ for a.e. $t\in[0,T]$. This fact, along with \emph{Step 1}, implies that $C_j\subset \mathscr{U}_1$, $\forall j\in\mathbb{N}$. \\
Since $X$ is a countable union of $\mathscr{Z}$-compact sets contained in $\mathscr{U}_1$ and $\mathscr{D}:\mathscr{Z}\rightarrow H^{-\varepsilon}\times\mathscr{V}'=:\mathscr{Y}$ is continuous as observed in Remark \eqref{rem:13}, it follows that $X\in\mathcal{B}_{\mathscr{Z}}(\mathscr{U}_1)$ and $\mathscr{F}:=\mathscr{D}(X)\in\mathcal{B}(\mathscr{Y})$. Moreover
\begin{equation}\label{eq:361}
(\Xi\times\mathbf{N})(\mathscr{F})=\mathbf{P}(\mathscr{D}^{-1}(\mathscr{F}))\geq\mathbf{P}(X)=1,
\end{equation}
where we have used \eqref{eq:330} from Lemma \eqref{lem:70}. Let $\mathscr{D}_1:X\rightarrow\mathscr{F}$ be the restriction $\mathscr{D}|_{X}$. Since $X\subset\mathscr{U}_1$, Corollary \eqref{cor:2} implies that $\mathscr{D}_1$ is one-to-one.
\vspace{1 pc}\\
\emph{Step 3: Measurability of $\mathscr{D}^{-1}_1$}. The mapping $\mathscr{D}^{-1}_1:(\mathscr{F},\mathcal{B}_{\mathscr{Y}}(\mathscr{F}))\rightarrow (X,\mathcal{B}_{\mathscr{Z}}(X))$ is measurable. 
To prove this, we only need to show that $\mathscr{D}(B)\in\mathcal{B_{\mathscr{Y}}}(\mathscr{F})$ for every $B$ closed set in $X$ in the topology of $\mathscr{Z}$, i.e. $B=B_1\cap X$ for $B_1$ closed set in $\mathscr{Z}$. Since $X=\cup_{j\in\mathbb{N}}{C_{j}}$, we have 
$$
\mathscr{D}(B)=\mathscr{D}(B_1\cap X)=\mathscr{D}(B_1\cap(\cup_{j\in\mathbb{N}}{C_{j}}))=\mathscr{D}(\cup_{j\in\mathbb{N}}{(B_1\cap C_{j}))}=\cup_{j\in\mathbb{N}}{\mathscr{D}(B_1\cap C_{j})}\in\mathcal{B}_{\mathscr{Y}}(\mathscr{F}),
$$
since $B_1\cap C_{j}$ is compact for each $j\in\mathbb{N}$, $\mathscr{D}:\mathscr{Z}\rightarrow\mathscr{Y}$ is continuous and hence $\mathscr{D}(B_1\cap C_{j})\in\mathcal{B}_{\mathscr{Y}}(\mathscr{F})$.
\vspace{1pc}\\
\emph{Step 4: Construction of the unique strong solution $\phi$}. We now set
$$
\Omega_1:=\left\{\omega\in\Omega:\left\{\phi_0(\omega),\frac{\partial w}{\partial t}(\omega)\right\}\in\mathscr{F}\right\}.
$$
Thanks to the measurability of $\phi_0$, $\partial w/\partial t$, and to Remark \eqref{rem:20}, we have $\Omega_1\in\mathcal{F}$. Relying on \eqref{eq:361}, we get 
$$
\mathbf{m}(\Omega_1)=(\Xi\times\mathbf{N})(\mathscr{F})=1.
$$
We have thus shown that $\mathscr{D}_1$ has the initial condition $\phi_0$ and the white noise $\partial w/\partial t$ in its image with full $\mathbf{m}$-probability, as anticipated in Section \ref{sec:41}. We finally define
$$
\phi(\omega):=
     \begin{cases}
      \mathscr{D}_1^{-1}\!\left(\left\{\phi_0(\omega),\displaystyle\frac{\partial w}{\partial t}(\omega)\right\}\right)\!, & \mbox{ if } \omega\in\Omega_1,\\
     0, & \mbox{ otherwise.}
     \end{cases}
$$
It is clear that $\phi$ satisfies \eqref{eq:331}. Since $X\in\mathcal{B}_{\mathscr{Z}}(\mathscr{U}_1)$, for each $G\in\mathcal{B}_{\mathscr{Z}}(\mathscr{U}_1)$ we have that $G\cap X\in\mathcal{B}_{\mathscr{Z}}(X)$. In addition, we get
\begin{eqnarray*}
\left\{\omega\in\Omega:\phi(\omega)\in G\right\} & = & \left\{\omega\in\Omega_1:\mathscr{D}^{-1}_1\left(\left\{\phi_0(\omega),\frac{\partial w}{\partial t}(\omega)\right\}\right)\in G\cap X\right\}\cup\left\{\omega\in\Omega^{C}_1:0\in G\right\}\\
& = &\underbrace{\left\{\omega\in\Omega_1:\left\{\phi_0(\omega),\frac{\partial w}{\partial t}(\omega)\right\}\in\mathscr{D}_1(G\cap X)\right\}}_{\displaystyle\Omega_2}\cup\underbrace{\left\{\omega\in\Omega^{C}_1:0\in G\right\}}_{\displaystyle\Omega_3}\in\mathcal{F}.
\end{eqnarray*} 
We have used the bijectivity of $\mathscr{D}^{-1}_{1}$ in the second equality, the measurability of $\mathscr{D}^{-1}_1$ to see that $\mathscr{D}_1(G\cap X)\in\mathcal{B}_{\mathscr{Y}}(\mathscr{F})$, the measurability of the random variable $\omega\mapsto\{\phi_0(\omega),\partial w/\partial t(\omega)\}$ to deduce that $\Omega_2\in\mathcal{F}$ and the fact that $\Omega_3$ is either $\emptyset$ or $\Omega^{C}_1$, hence $\Omega_3\in\mathcal{F}$ in both cases.\\
Hence, $\phi$ is measurable as a random variable from $(\Omega,\mathcal{F})$ to $(\mathscr{U}_1,\mathcal{B}_{\mathscr{Z}}(\mathscr{U}_1))$ and is then a strong solution in the sense of Definition \eqref{def:4}. In addition, the injectivity of $\mathscr{D}:\mathscr{U}_1\rightarrow\mathscr{Y}$ also implies the uniqueness of a strong solution.
\vspace{1pc}\\
\emph{Step 5: Distribution of $\phi$}. It is also clear, from the computations made in Remark \eqref{rem:16}, that the distribution of a strong solution is a weak solution. In addition, $\mathbf{P}$ is the distribution of $\phi$: to prove this, let $\mathbf{P}_1$ be another weak solution concentrated on $\mathcal{B}_{\mathscr{Z}}(X)$
. We rely on the bijectivity of $\mathscr{D}_1$ and write, for any given subset $C\in\mathcal{B}_{\mathscr{Z}}(X)$,
\begin{eqnarray*}
\mathbf{P}_1(C) & = & \mathbf{P}_1(\mathscr{D}^{-1}_1(\mathscr{D}_1(C)))=\mathbf{P}_1(\mathscr{D}^{-1}(\mathscr{D}_1(C)))=(\Xi\times\mathbf{N})(\mathscr{D}_1(C))\\
& = & \mathbf{P}(\mathscr{D}^{-1}(\mathscr{D}_1(C)))=\mathbf{P}(\mathscr{D}^{-1}_1(\mathscr{D}_1(C)))=\mathbf{P}(C).
\end{eqnarray*}
Hence $\mathbf{P}=\mathbf{P}_1$. Since the distribution of $\phi$ is clearly a weak solution concentrated on $\mathcal{B}_{\mathscr{Z}}(X)$, we deduce that $\mathbf{P}$ is the distribution of $\phi$. 
\end{Proof}
\begin{rem}
We point out that \emph{Proof of Theorem} \eqref{th:18}, \emph{Step 5} can be seen as a \emph{partial} uniqueness result for weak solutions, since the requirement upon the $\sigma$-algebra $\mathcal{B}_{\mathscr{Z}}(X)$ restricts the set of probability measures in which we look for a weak solution. 
\end{rem}
\begin{rem}
The spatial dimension requirement $d\leq 3$ is needed only in \emph{Proof of Theorem} \eqref{th:20} to guarantee the compact embedding $U\hookrightarrow L^4$ and in \emph{Proof of Theorem} \eqref{th:18}, \emph{Step 2}, to compute the distributional gradient of $v^3$, for any $v\in U$.
\end{rem}
\begin{rem}\label{rem:199}
Apparently it is not straightforward to relax the polynomial growth restriction on $F$
 we have imposed in Section \ref{sec:5} in order to be able to prove Theorems \eqref{th:19} and \eqref{th:18}. 
\end{rem}
\thanks{{\scshape Acknowledgements}. This paper is based on the author's Master thesis defended at the Politecnico di Milano (Italy) under the supervision of Maurizio Grasselli, with the co-supervision of Marco Alessandro Fuhrman. The author's deep gratitude goes primarily to both of them for suggesting the problem and for their constant guidance and encouragement. The author also wishes to thank Giuseppe Da Prato and Benedetta Ferrario for their help on some specific points.}

\bibliographystyle{amsplain}

{\footnotesize
\bibliography{Biblio}}

\providecommand{\bysame}{\leavevmode\hbox to3em{\hrulefill}\thinspace}
\providecommand{\MR}{\relax\ifhmode\unskip\space\fi MR }
\providecommand{\MRhref}[2]{%
  \href{http://www.ams.org/mathscinet-getitem?mr=#1}{#2}
}
\providecommand{\href}[2]{#2}
\begin{thebibliography}{10}

\bibitem{AK}
D.~Antonopoulou and G.~Karali, \emph{Existence of solution for a generalized
  stochastic {C}ahn-{H}illiard equation on convex domains}, Discrete Contin.
  Dyn. Syst. Ser. B \textbf{16} (2011), no.~1, 31--55.

\bibitem{AG}
A.~Araujo and E.~Gin\'e, \emph{The central limit theorem for real and {B}anach
  valued random variables}, vol. 431, Wiley, New York, 1980.

\bibitem{BH2}
P.~W. Bates and J.{} Han, \emph{The {D}irichlet boundary problem for a nonlocal
  {C}ahn-{H}illiard equation}, J. Math. Anal. Appl. \textbf{311} (2005), no.~1,
  289--312.

\bibitem{BH1}
P.~W. Bates and J.~Han, \emph{The {N}eumann boundary problem for a nonlocal
  {C}ahn-{H}illiard equation}, J. Differential Equations \textbf{212} (2005),
  no.~2, 235--277.

\bibitem{CH1}
J.~W. Cahn and J.~E. Hilliard, \emph{Free energy of a nonuniform system {I}.
  interfacial free energy}, J. Chem. Phys. \textbf{28} (1958), no.~2, 258--267.

\bibitem{CW}
C.~Cardon-Weber, \emph{Cahn-{H}illiard stochastic equation: existence of the
  solution and of its density}, Bernoulli \textbf{7} (2001), no.~5, 777--816.

\bibitem{CFG}
P.~Colli, S.~Frigeri, and M.~Grasselli, \emph{Global existence of weak
  solutions to a nonlocal {C}ahn-{H}illiard-{N}avier-{S}tokes system}, J. Math.
  Anal. Appl. \textbf{386} (2012), no.~1, 428--444.

\bibitem{CO}
H.~E. Cook, \emph{Brownian motion in spinodal decomposition}, Acta Metallurgica
  \textbf{18} (1970), no.~3, 297--306.

\bibitem{DPD}
G.~Da~Prato and A.~Debussche, \emph{Stochastic {C}ahn-{H}illiard equation},
  Nonlinear Anal. \textbf{26} (1996), no.~2, 241--263.

\bibitem{DPZ}
G.~Da~Prato and J.~Zabczyk, \emph{Stochastic equations in infinite dimensions},
  Encyclopedia of Mathematics and its Applications, Cambridge University Press,
  Cambridge, 1992.

\bibitem{DZ}
A.~Debussche and L.~Zambotti, \emph{Conservative stochastic {C}ahn-{H}illiard
  equation with reflection}, Ann. Probab. \textbf{35} (2007), no.~5,
  1706--1739.

\bibitem{LHJSYK}
L.~Dongsun, H.~Joo-Youl, J.~Darae, S.~Jaemin, Y.~Ana, and K.~Junseok,
  \emph{Physical, mathematical, and numerical derivations of the
  {C}ahn-{H}illiard equation}, Computational Materials Science \textbf{81}
  (2014), 216--225.

\bibitem{EM}
N.~Elezovi\'c and A.~Mikeli\'c, \emph{On the {S}tochastic {C}ahn-{H}illiard
  {E}quation}, Nonlinear Anal. \textbf{16} (1991), no.~12, 1169--1200.

\bibitem{FG}
S.~Frigeri and M.~Grasselli, \emph{Global and {T}rajectory {A}ttractors for a
  {N}onlocal {C}ahn-{H}illiard-{N}avier-{S}tokes {S}ystem}, J. Dynam.
  Differential Equations \textbf{24} (2012), no.~4, 827--856.

\bibitem{FV}
A.~V. Fursikov and M.~I. Vishik, \emph{Mathematical {P}roblems of {S}tatistical
  {H}ydromechanics}, vol.~9, Kluver Academic Publishers, Dordrecht, 1988.

\bibitem{GAJ}
H.~Gajewski, \emph{On a nonlocal model of non-isothermal phase separation},
  Adv. Math. Sci. Appl. \textbf{12} (2002), no.~2, 569--586.

\bibitem{GZ}
H.~Gajewski and K.~Zacharias, \emph{On a nonlocal phase separation model}, J.
  Math. Anal. Appl. \textbf{286} (2003), no.~1, 11--31.

\bibitem{GL1}
G.~Giacomin and J.~L. Lebowitz, \emph{Phase segregation dynamics in particle
  systems with long range interactions. {I}. {M}acroscopic limits}, J. Statist.
  Phys. \textbf{87} (1997), no.~1-2, 37--61.

\bibitem{GL2}
G.~{G}iacomin and J.~L. Lebowitz, \emph{Phase segregation dynamics in particle
  systems with long range interactions. {II}. {I}nterface motion}, SIAM J.
  Appl. Math. \textbf{58} (1998), no.~6, 1707--1729.

\bibitem{GR}
D.~Grieser, \emph{Uniform bounds for eigenfunctions of the {L}aplacian on
  manifolds with boundary}, Comm. Partial Differential Equations \textbf{27}
  (2002), no.~7-8, 1283--1299.

\bibitem{HIL}
J.~E. {H}illiard, \emph{Spinodal decomposition, \emph{paper from {P}hase
  {T}ransformations}}, ASM (1970), 497--560.

\bibitem{IW}
N.~Ikeda and Sh. Watanabe, \emph{Stochastic {D}ifferential {E}quations and
  {D}iffusion {P}rocesses}, North-Holland, Amsterdam, 1981.

\bibitem{JP}
J.~Jacod and P.~Protter, \emph{Probability {E}ssentials}, second ed.,
  Springer-Verlag, Berlin, 2004, (corrected second printing).

\bibitem{LM}
J.-L. Lions and E.~Magenes, \emph{Non-homogeneous boundary value problem and
  applications}, vol. I and II, Grundlehren der mathematischen Wissenschaften,
  no. 181, Springer-Verlag, Berlin, 1972.

\bibitem{PR}
C.~Pr\'evot and M.~R\"ockner, \emph{A {C}oncise {C}ourse on {S}tochastic
  {P}artial {D}ifferential {E}quations}, Springer-Verlag, Berlin, 2007.

\bibitem{RW}
L.~C.~G. Rogers and D.~Williams, \emph{Diffusions, {M}arkov {P}rocesses and
  {M}artingales}, second ed., Cambridge University Press, Cambridge, 2000,
  (reissued).

\bibitem{S}
J.~Simon, \emph{Compact sets in the space {$L^{p}(0,T;B)$}}, Annali Mat. Pura
  Appl. \textbf{CXLVI} (1987), no.~1, 65--96.

\end{thebibliography}


%
%
%
%
%



\end{document}